\documentclass[a4paper,12pt]{article}
\usepackage{graphicx}
\usepackage{amsthm,amsmath, amsfonts}
\usepackage{amsthm,amscd,amsfonts,latexsym,color,epsfig}
\usepackage{amssymb, epic, graphicx}

\begin{document}

\sloppy

\title{Inverse source non-local problem for mixed type equation with Caputo fractional differential operator}
\author{Erkinjon Karimov, Nasser Al-Salti,  and Sebti Kerbal}

\maketitle

\begin{abstract} 

In the present work, we discuss a unique solvability of an inverse-source problem with integral transmitting condition for time-fractional mixed type equation in a rectangular domain, where the unknown source term depends on space variable only.

The method of solution based on a series expansion using bi-orthogonal basis of space corresponding to a nonself-adjoint boundary value problem. Under certain regularity conditions on the given data, we prove the uniqueness and existence of the solution for the given problem. Influence of transmitting condition on the solvability of the problem is shown as well. Precisely, two different cases were considered; a case of full integral form ($0<\gamma<1$) and a special case ($\gamma=1$) of transmitting condition.

In order to simplify the bulky expressions appearing in the proof of the main result, we have  established a new property of the recently introduced Mittag-Leffler type function of two variables (see Lemma 2.1).

\end{abstract}



\section{Introduction}

Fractional differential equations (FDEs) have become an important object in modeling many real life processes arising in different fields such as water movement in soil [1], nanotechnology [2] and viscoelasticity [3]. Especially, inverse problems (IPs) involving FDEs have numerous applications. It is known that problems of identifying coefficients or order of considered equation, boundary conditions and/or source function are considered as IPs. Latter case is known as inverse source problem, where one is looking for the source term(s) in a differential equation (or a system of equations) using extra boundary data. For example, in [4], authors dealt with an inverse problem of simultaneously identifying the space-dependent diffusion coefficient and the fractional order  in the 1D time-fractional diffusion equation. Kirane et al [5] studied conditional well posedness in determining a space-dependent source in the 2D time-fractional diffusion equation. In [6], the source term for a time fractional diffusion equation was determined with an integral type over-determining condition.

Inverse problems for various mixed type equations of integer order were studied by Sabitov and his students in several works. For instance, in [7] they have considered the following Lavrent'ev-Bitsadze equation
$$
{{u}_{xx}}+sign(y){{u}_{yy}}=f\left( x \right)
$$
in a rectangular domain. Using the method of spectral expansions, they proved the unique solvability of inverse problem with non-local conditions with respect to space variables. Another problem on source identification was investigated in [8], precisely, well-posedness and stability inequalities for the solution of three source identification problems for hyperbolic-parabolic equations were obtained.

In [9], inverse problems for time-fractional mixed type equations with uniformly elliptic operator with respect to space variable were studied in a rectangular domain. Under certain regularity conditions to the given functions and geometric restrictions to the considered domain, a unique weak solvability was proved.

Influence of transmitting conditions in studying mixed type equations becomes interesting, since depending on the form of transmitting conditions, restrictions on the given data could be reduced. For instance, see [10] for fractional case and [11] for integer case of mixed type equations.

In the present work, we investigate inverse source problem with non-local conditions for time-fractional mixed type equation with Caputo derivative. Due to non-local condition, we have used the method of series expansion using a bi-orthogonal basis corresponding to a nonself-adjoint boundary value problem. We note that this kind of bi-orthogonal basis was used in [12], where the studied equation contains generalized fractional derivative.

The rest of the paper is organized as follows. In section 2, we give preliminary information on definition, properties of Euler's integrals, Mittag-Leffler type functions, Riemann-Liouville, Caputo fractional differential operators and a short info on the used bi-orthogonal system as well. Section 3 is devoted to the formulation of the problem and construction of formal solution. In Section 4, we have proved the uniqueness and the existence of solution for the case $0<\gamma<1$. In Section 5, special case of transmitting condition $\gamma=1$ is considered. At the end of paper, some details of calculations, used in the paper are presented in the Appendix.

\section{Preliminaries}
\subsection{Euler's integrals}

The elementary definition of the gamma function is the Euler's integral [13, p.24]:
$$
\Gamma \left( z \right)=\int\limits_{0}^{\infty }{{{t}^{z-1}}{{e}^{-t}}dt}.
$$
For $z\in {{\mathbb{R}}^{+}}$, this integral converges and satisfies the recurrence relation
$$
\Gamma \left( z+1 \right)=z\Gamma \left( z \right)\eqno (2.1)
$$
and for $z\in \mathbb{N}$ it satisfies
$$
\Gamma \left( 1+z \right)=z! \eqno (2.2)
$$

There is another Euler's integral, which can be represented by gamma function [13, p.26]:
$$
\frac{\Gamma \left( x \right)\Gamma \left( y \right)}{\Gamma \left( x+y \right)}=\int\limits_{0}^{1}{{{t}^{x-1}}{{\left( 1-t \right)}^{y-1}}dt}, \eqno (2.3)
$$
which converges for $x>0,\,y>0$.

\subsection{Mittag-Leffler type functions}

A two-parameter function of the Mittag-Leffler type is defined by the series expansion [14, p.17] as follows
$$
{{E}_{\alpha ,\beta }}\left( z \right)=\sum\limits_{k=0}^{\infty }{\frac{{{z}^{k}}}{\Gamma \left( \alpha k+\beta  \right)}}\,\,\left( \alpha >0,\,\beta >0 \right),\eqno (2.4)
$$
which satisfies the following relation [13,p.45] and formula of differentiation [14, p.21], respectively
$$
{{E}_{\alpha ,\beta }}\left( z \right)-z{{E}_{\alpha ,\alpha +\beta }}\left( z \right)=\frac{1}{\Gamma \left( \beta  \right)},\eqno (2.5)
$$
and
$$
{{\,}_{RL}}D_{0t}^{\gamma }\left( {{t}^{\alpha k+\beta -1}}E_{\alpha ,\beta }^{\left( k \right)}\left( \lambda {{t}^{\alpha }} \right) \right)={{t}^{\alpha k+\beta -\gamma -1}}E_{\alpha ,\beta -\gamma }^{\left( k \right)}\left( \lambda {{t}^{\alpha }} \right), \eqno (2.6)
$$
where $E_{\alpha ,\beta }^{(k)}(t)=\frac{{{d}^{k}}}{d{{t}^{k}}}{{E}_{\alpha ,\beta }}(t)$, denotes the classical derivative of order $k$.

It also satisfies the inequality presented in the following theorem.

\textbf{Theorem 2.1.} (Theorem 1.6 in [14]) If $\alpha <2$, $\beta $ is an arbitrary real number, $\mu $ is a real number such that $\pi \alpha /2<\mu <\min \{\pi ,\pi \alpha \}$ and $C$ is a real constant, then
$$
\left| {{E}_{\alpha ,\beta }}\left( z \right) \right|\le \frac{C}{1+\left| z \right|}, \,\,\,\,\left( \mu \le \arg z\le \pi  \right),\,\,\left| z \right|\ge 0.
$$

The following Mittag-Leffler type function in two variables
$$
\begin{aligned}{l}
&{{E}_{1}}\left( \begin{matrix}
   {{\gamma }_{1}},{{\alpha }_{1}};{{\gamma }_{2}},{{\beta }_{1}} & |\,\,x  \\
   {{\delta }_{1}},{{\alpha }_{2}},{{\beta }_{2}};{{\delta }_{2}},{{\alpha }_{3}};{{\delta }_{3}},{{\beta }_{3}} & |\,\,y  \\
\end{matrix} \right)=\sum\limits_{m=0}^{\infty }{\sum\limits_{n=0}^{\infty }{\frac{{{\left( {{\gamma }_{1}} \right)}_{{{\alpha }_{1}}m}}{{\left( {{\gamma }_{2}} \right)}_{{{\beta }_{1}}n}}}{\Gamma \left( {{\delta }_{1}}+{{\alpha }_{2}}m+{{\beta }_{2}}n \right)}\times}}\\
&\times \frac{{{x}^{m}}}{\Gamma \left( {{\delta }_{2}}+{{\alpha }_{3}}m \right)}\frac{{{y}^{n}}}{\Gamma \left( {{\delta }_{3}}+{{\beta }_{3}}n \right)}\\
\end{aligned}
\eqno (2.7)
$$
was introduced by Garg et al [15] in 2013, where $(\gamma)_n=\frac{\Gamma(\gamma+n)}{\Gamma(\gamma)}$ is Pochgammer's symbol. This function has the following integral representation
$$
\begin{aligned}
  & {{E}_{1}}\left( \begin{matrix}
   {{\gamma }_{1}},{{\alpha }_{1}};{{\gamma }_{2}},{{\beta }_{1}} & |\,\,x  \\
   {{\delta }_{1}},{{\alpha }_{2}},{{\beta }_{2}};{{\delta }_{2}},{{\alpha }_{3}};{{\delta }_{3}},{{\beta }_{3}} & |\,\,y  \\
\end{matrix} \right)= \frac{1}{\Gamma \left( {{\gamma }_{1}} \right)\Gamma \left( {{\gamma }_{2}} \right)}\times\\
 & \times\int\limits_{0}^{1}{{{t}^{{{\rho }_{1}}-1}}{{\left( 1-t \right)}^{{{\rho }_{2}}-1}}E_{{{\alpha }_{2}},{{\rho }_{1}};{{\alpha }_{3}},{{\delta }_{2}}}^{{{\gamma }_{1}},{{\alpha }_{1}}}\left( x{{t}^{{{\alpha }_{2}}}} \right)E_{{{\beta }_{2}},{{\rho }_{2}};{{\beta }_{3}},{{\delta }_{3}}}^{{{\gamma }_{2}},{{\beta }_{1}}}\left( y{{\left( 1-t \right)}^{{{\beta }_{2}}}} \right)dt} \\
\end{aligned} \eqno (2.8)
$$
and the following formula for differentiation [15]
$$
\begin{aligned}
  & _{RL}D_{yo}^{\gamma }\left\{ {{\left( -y \right)}^{{{\delta }_{1}}-1}}{{E}_{1}}\left( \begin{matrix}
   {{\gamma }_{1}},{{\alpha }_{1}};{{\gamma }_{2}},{{\beta }_{1}} & |\,\,{{\omega }_{1}}{{\left( -y \right)}^{{{\alpha }_{2}}}}  \\
   {{\delta }_{1}},{{\alpha }_{2}},{{\beta }_{2}};{{\delta }_{2}},{{\alpha }_{3}};{{\delta }_{3}},{{\beta }_{3}} & |\,\,{{\omega }_{2}}{{\left( -y \right)}^{{{\beta }_{2}}}}  \\
\end{matrix} \right) \right\}= \\
 & ={{\left( -y \right)}^{{{\delta }_{1}}-\gamma -1}}{{E}_{1}}\left( \begin{matrix}
   {{\gamma }_{1}},{{\alpha }_{1}};{{\gamma }_{2}},{{\beta }_{1}} & |\,\,{{\omega }_{1}}{{\left( -y \right)}^{{{\alpha }_{2}}}}  \\
   {{\delta }_{1}}-\gamma ,{{\alpha }_{2}},{{\beta }_{2}};{{\delta }_{2}},{{\alpha }_{3}};{{\delta }_{3}},{{\beta }_{3}} & |\,\,{{\omega }_{2}}{{\left( -y \right)}^{{{\beta }_{2}}}}  \\
\end{matrix} \right). \\
\end{aligned}\eqno (2.9)
$$
Here
$$
E_{{{\alpha }_{2}},{{\delta }_{1}};{{\alpha }_{3}},{{\delta }_{2}}}^{{{\gamma }_{1}},{{\alpha }_{1}}}\left( x \right)=\sum\limits_{m=0}^{\infty }{\frac{{{\left( {{\gamma }_{1}} \right)}_{{{\alpha }_{1}}m}}}{\Gamma \left( {{\delta }_{1}}+{{\alpha }_{2}}m \right)}\frac{{{x}^{m}}}{\Gamma \left( {{\delta }_{2}}+{{\alpha }_{3}}m \right)}}\eqno (2.10)
$$
is the Mittag-Leffler type function in one variable and
$$
\begin{aligned}
&{{\gamma }_{1}},\,{{\gamma }_{2}},\,{{\delta }_{1}},\,{{\delta }_{2}},\,{{\delta }_{3}},\,x,\,y\in \mathbb{C},\,\min \left\{ {{\alpha }_{1}},\,{{\alpha }_{2}},\,{{\alpha }_{3}},\,\,{{\beta }_{1}},\,{{\beta }_{2}},\,{{\beta }_{3}},\right.\\
&\left.\operatorname{Re}\left( {{\rho }_{1}} \right),\,\operatorname{Re}\left( {{\rho }_{2}} \right) \right\}>0,\,{{\rho }_{1}}+{{\rho }_{2}}={{\delta }_{1}}.\\
\end{aligned}
$$

We have also found that this function satisfies a useful relation, which is presented in the following lemma.

\textbf{Lemma 2.2.} If $y=x,\,{{\gamma }_{2}}={{\delta }_{3}},\,{{\beta }_{1}}={{\beta }_{3}},\,\lambda ={{\beta }_{2}}$, then
$$
\begin{aligned}
&{{E}_{1}}\left( \begin{matrix}
   {{\gamma }_{1}},{{\alpha }_{1}};{{\gamma }_{2}},{{\beta }_{1}} & |\,\,x  \\
   {{\delta }_{1}},{{\alpha }_{2}},{{\beta }_{2}};{{\delta }_{2}},{{\alpha }_{3}};{{\delta }_{3}},{{\beta }_{3}} & |\,\,y  \\
\end{matrix} \right)-y{{E}_{1}}\left( \begin{matrix}
   {{\gamma }_{1}},{{\alpha }_{1}};{{\gamma }_{2}},{{\beta }_{1}} & |\,\,x  \\
   {{\delta }_{1}}+\lambda ,{{\alpha }_{2}},{{\beta }_{2}};{{\delta }_{2}},{{\alpha }_{3}};{{\delta }_{3}},{{\beta }_{3}} & |\,\,y  \\
\end{matrix} \right)=\\
&=E_{{{\alpha }_{2}},{{\delta }_{1}};{{\alpha }_{3}},{{\delta }_{2}}}^{{{\gamma }_{1}},{{\alpha }_{1}}}\left( x \right).\\
\end{aligned}
\eqno (2.11)
$$

\texttt{Proof:}

Using representation (2.7), we get
$$
\begin{aligned}
&\sum\limits_{m=0}^{\infty }{\sum\limits_{n=0}^{\infty }{\frac{{{\left( {{\gamma }_{1}} \right)}_{{{\alpha }_{1}}m}}{{\left( {{\gamma }_{2}} \right)}_{{{\beta }_{1}}n}}}{\Gamma \left( {{\delta }_{1}}+{{\alpha }_{2}}m+{{\beta }_{2}}n \right)}\frac{{{x}^{m}}}{\Gamma \left( {{\delta }_{2}}+{{\alpha }_{3}}m \right)}\frac{{{y}^{n}}}{\Gamma \left( {{\delta }_{3}}+{{\beta }_{3}}n \right)}}}-\\
&-x\sum\limits_{m=0}^{\infty }{\sum\limits_{n=0}^{\infty }{\frac{{{\left( {{\gamma }_{1}} \right)}_{{{\alpha }_{1}}m}}{{\left( {{\gamma }_{2}} \right)}_{{{\beta }_{1}}n}}}{\Gamma \left( {{\delta }_{1}}+\lambda +{{\alpha }_{2}}m+{{\beta }_{2}}n \right)}\frac{{{x}^{m}}}{\Gamma \left( {{\delta }_{2}}+{{\alpha }_{3}}m \right)}\frac{{{y}^{n}}}{\Gamma \left( {{\delta }_{3}}+{{\beta }_{3}}n \right)}}}. \\
\end{aligned}
$$
On setting $y=x$, the above equation can be written as
$$
\sum\limits_{m=0}^{\infty }{\frac{{{\left( {{\gamma }_{1}} \right)}_{{{\alpha }_{1}}m}}{{x}^{m}}}{\Gamma \left( {{\delta }_{2}}+{{\alpha }_{3}}m \right)}\left[ \frac{{{\left( {{\gamma }_{2}} \right)}_{0}}}{\Gamma \left( {{\delta }_{1}}+{{\alpha }_{2}}m \right)\Gamma \left( {{\delta }_{3}} \right)}+\frac{{{\left( {{\gamma }_{2}} \right)}_{{{\beta }_{1}}}}x}{\Gamma \left( {{\delta }_{1}}+{{\alpha }_{2}}m+{{\beta }_{2}} \right)\Gamma \left( {{\delta }_{3}}+{{\beta }_{3}} \right)} \right.+}
$$
$$
\begin{aligned}
  & +\frac{{{\left( {{\gamma }_{2}} \right)}_{2{{\beta }_{1}}}}{{x}^{2}}}{\Gamma \left( {{\delta }_{1}}+{{\alpha }_{2}}m+2{{\beta }_{2}} \right)\Gamma \left( {{\delta }_{3}}+2{{\beta }_{3}} \right)}+\frac{{{\left( {{\gamma }_{2}} \right)}_{3{{\beta }_{1}}}}{{x}^{3}}}{\Gamma \left( {{\delta }_{1}}+{{\alpha }_{2}}m+3{{\beta }_{2}} \right)\Gamma \left( {{\delta }_{3}}+3{{\beta }_{3}} \right)}+...\\
  &-\frac{{{\left( {{\gamma }_{2}} \right)}_{0}}x}{\Gamma \left( {{\delta }_{1}}+\lambda +{{\alpha }_{2}}m \right)\Gamma \left( {{\delta }_{3}} \right)}-\frac{{{\left( {{\gamma }_{2}} \right)}_{{{\beta }_{1}}}}{{x}^{2}}}{\Gamma \left( {{\delta }_{1}}+\lambda +{{\alpha }_{2}}m+{{\beta }_{2}} \right)\Gamma \left( {{\delta }_{3}}+{{\beta }_{3}} \right)}- \\
 & \left. -\frac{{{\left( {{\gamma }_{2}} \right)}_{2{{\beta }_{1}}}}{{x}^{3}}}{\Gamma \left( {{\delta }_{1}}+\lambda +{{\alpha }_{2}}m+2{{\beta }_{2}} \right)\Gamma \left( {{\delta }_{3}}+2{{\beta }_{3}} \right)}-\frac{{{\left( {{\gamma }_{2}} \right)}_{3{{\beta }_{1}}}}{{x}^{4}}}{\Gamma \left( {{\delta }_{1}}+\lambda +{{\alpha }_{2}}m+3{{\beta }_{2}} \right)\Gamma \left( {{\delta }_{3}}+3{{\beta }_{3}} \right)}-... \right],\\
\end{aligned}
$$
which on assuming ${{\gamma }_{2}}={{\delta }_{3}},\,{{\beta }_{1}}={{\beta }_{3}},\lambda ={{\beta }_{2}}$ becomes
$$
\sum\limits_{m=0}^{\infty }{\frac{{{\left( {{\gamma }_{1}} \right)}_{{{\alpha }_{1}}m}}{{x}^{m}}}{\Gamma \left( {{\delta }_{2}}+{{\alpha }_{3}}m \right)\Gamma \left( {{\delta }_{1}}+{{\alpha }_{2}}m \right)}=}E_{{{\alpha }_{2}},{{\delta }_{1}};{{\alpha }_{3}},{{\delta }_{2}}}^{{{\gamma }_{1}},{{\alpha }_{1}}}\left( x \right).
$$
As a special case, we obtain
$$
{{E}_{1}}\left( \begin{matrix}
   1,1;1,1 & |\,\,x  \\
   \alpha +1,\alpha ,\alpha ;1,1;1,1 & |\,\,x  \\
\end{matrix} \right)-x{{E}_{1}}\left( \begin{matrix}
   1,1;1,1 & |\,\,x  \\
   \alpha +1,\alpha ,\alpha ;1,1;1,1 & |\,\,x  \\
\end{matrix} \right)={{E}_{\alpha ,\alpha +1}}\left( x \right), \eqno (2.12)
$$
which will be used later in simplifying our calculations. 

\subsection{Riemann-Liouville and Caputo fractional derivatives}

The Riemann-Liouville fractional derivatives $_{RL}D_{ax}^{\alpha }\,y$ and $_{RL}D_{xb}^{\alpha }\,y$ of order $\alpha \in \mathbb{C}\,\left( \operatorname{Re}\left( \alpha  \right)\ge 0 \right)$ are defined by [13, p.70]
$$
\begin{aligned}
  & \left( _{RL}D_{ax}^{\alpha }\,y \right)\left( x \right)=\frac{1}{\Gamma \left( n-\alpha  \right)}{{\left( \frac{d}{dx} \right)}^{n}}\int\limits_{a}^{x}{\frac{y\left( t \right)dt}{{{\left( x-t \right)}^{\alpha -n+1}}}\,\,\,\,\,\,\,\,\,\,\left( n=\left[ \operatorname{Re}\left( \alpha  \right) \right]+1,\,x>a \right)}, \\
 & \left( _{RL}D_{xb}^{\alpha }\,y \right)\left( x \right)=\frac{1}{\Gamma \left( n-\alpha  \right)}{{\left( -\frac{d}{dx} \right)}^{n}}\int\limits_{x}^{b}{\frac{y\left( t \right)dt}{{{\left( t-x \right)}^{\alpha -n+1}}}\,\,\,\,\,\,\,\left( n=\left[ \operatorname{Re}\left( \alpha  \right) \right]+1,\,x<b \right)}, \\
\end{aligned} \eqno (2.13)
$$
respectively, where $\left[ \operatorname{Re}\left( \alpha  \right) \right]$ is the integer part of $\operatorname{Re}\left( \alpha  \right)$. In particular, for $\alpha =n\in \mathbb{N}\cup \left\{ 0 \right\}$, we have
$$
\begin{aligned}
&\left( _{RL}D_{ax}^{0}\,y \right)\left( x \right)=y\left( x \right),\,\left( _{RL}D_{xb}^{0}\,y \right)\left( x \right)=y\left( x \right),\,\left( _{RL}D_{ax}^{n}\,y \right)\left( x \right)={{y}^{\left( n \right)}}\left( x \right),\\
&\left( _{RL}D_{xb}^{n}\,y \right)\left( x \right)={{\left( -1 \right)}^{n}}{{y}^{\left( n \right)}}\left( x \right),\,n\in \mathbb{N},
\end{aligned}
$$
where ${{y}^{\left( n \right)}}\left( x \right)$ is the usual derivative of $y\left( x \right)$ of order $n$.

If $\alpha \notin \mathbb{N}\cup \left\{ 0 \right\}$, the Caputo fractional derivatives $_{C}D_{ax}^{\alpha }\,y$ and $_{C}D_{xb}^{\alpha }\,y$ of order $\alpha$   are defined by [13,p.92]
$$
\begin{aligned}
  & \left( _{C}D_{ax}^{\alpha }\,y \right)\left( x \right)=\frac{1}{\Gamma \left( n-\alpha  \right)}\int\limits_{a}^{x}{\frac{{{y}^{\left( n \right)}}\left( t \right)dt}{{{\left( x-t \right)}^{\alpha -n+1}}}}\,\,\,\,\,\,\,\,\,\,\,\,\left( n=\left[ \operatorname{Re}\left( \alpha  \right) \right]+1,\,x>a \right),\, \\
 & \left( _{C}D_{xb}^{\alpha }\,y \right)\left( x \right)=\frac{{{\left( -1 \right)}^{n}}}{\Gamma \left( n-\alpha  \right)}\int\limits_{x}^{b}{\frac{{{y}^{\left( n \right)}}\left( t \right)dt}{{{\left( t-x \right)}^{\alpha -n+1}}}\,}\,\,\,\,\,\,\,\,\,\,\,\left( n=\left[ \operatorname{Re}\left( \alpha  \right) \right]+1,\,x<b \right), \\
\end{aligned}\eqno (2.14)
$$
respectively, while for $\alpha =n\in \mathbb{N}\cup \left\{ 0 \right\}$, we have
$$
\begin{aligned}
& \left( _{C}D_{ax}^{0}y \right)\left( x \right)=y\left( x \right),\,\left( _{C}D_{xb}^{0}y \right)\left( x \right)=y\left( x \right),\,\left( _{C}D_{ax}^{n}y \right)\left( x \right)={{y}^{\left( n \right)}}\left( x \right),\\
& \,\left( _{C}D_{xb}^{n}y \right)\left( x \right)={{\left( -1 \right)}^{n}}{{y}^{\left( n \right)}}\left( x \right).\\
\end{aligned}
$$

These two operators are connected with each other by the following relations [13,p.91]
$$
\begin{aligned}
  & \left( _{C}D_{ax}^{\alpha }\,y \right)\left( x \right)=\left( _{RL}D_{ax}^{\alpha }\,y \right)\left( x \right)-\sum\limits_{k=0}^{n-1}{\frac{{{y}^{\left( k \right)}}\left( a \right)}{\Gamma \left( k-\alpha +1 \right)}{{\left( x-a \right)}^{k-\alpha }}}\\
 & \left( n=\left[ \operatorname{Re}\left( \alpha  \right) \right]+1,\,x>a \right), \\
 & \left( _{C}D_{xb}^{\alpha }\,y \right)\left( x \right)=\left( _{RL}D_{xb}^{\alpha }\,y \right)\left( x \right)-\sum\limits_{k=0}^{n-1}{\frac{{{y}^{\left( k \right)}}\left( b \right)}{\Gamma \left( k-\alpha +1 \right)}{{\left( b-x \right)}^{k-\alpha }}}\\
 & \left( n=\left[ \operatorname{Re}\left( \alpha  \right) \right]+1,\,x<b \right). \\
\end{aligned}\eqno (2.15)
$$

\subsection{Bi-orthogonal system}

Consider the following spectral problem:
$$
{X}''\left( x \right)+\mu X\left( x \right)=0,\,\,X\left( 0 \right)=X\left( 1 \right),\,\,{X}'\left( 0 \right)=0,\eqno (2.16)
$$
which has eigenvalues ${{\mu }_{k}}=\lambda _{k}^{2},\,\,{{\lambda }_{0}}=0,\,{{\lambda }_{k}}=2k\pi \,(k=1,2,...)$ and corresponding eigenfunctions $1,\,\cos {{\lambda }_{k}}x$, supplemented by the associate function $x\sin {{\lambda }_{k}}x$. The system of root functions is given by
$$
{{X}_{k}}\left( x \right)=\left\{ 1,\,\cos {{\lambda }_{k}}x,\,x\sin {{\lambda }_{k}}x \right\}.\eqno (2.17)
$$

Since, problem (2.16) is not self-adjoint, we need to find root functions of the corresponding adjoint problem, which is
$$
{Y}''\left( x \right)+\mu Y\left( x \right)=0,\,\,{Y}'\left( 0 \right)={Y}'\left( 1 \right),\,\,Y\left( 1 \right)=0.
$$
This problem has the same eigenvalues as problem (2.16), but another system of root functions, namely
$$
{{Y}_{k}}\left( x \right)=\left\{ 2(1-x),4(1-x)\,\cos {{\lambda }_{k}}x,\,4\sin {{\lambda }_{k}}x \right\}.\eqno (2.18)
$$

Systems (2.17)-(2.18) form bi-orthogonal system, which satisfies the necessary and sufficient condition for the basis property in the space ${{L}_{2}}[0,1]$ (see [16]).

\section{Problem formulation and formal representation of the solution}

Consider a time-fractional mixed type equation
$$
f\left( x \right)=\left\{ \begin{aligned}
  & _{C}D_{0t}^{\alpha }u-{{u}_{xx}},\,\,t>0, \\
 & _{C}D_{t0}^{\beta }u-{{u}_{xx}},\,\,t<0 \\
\end{aligned} \right.\eqno (3.1)
$$
in a rectangular domain $\Omega =\left\{ \left( x,t \right):\,0<x<1,\,-p<t<q \right\}.$  Here $\alpha ,\,\beta ,\,p,\,q\in {{\mathbb{R}}^{+}}$ such that $0<\alpha \le 1,\,1<\beta \le 2$, $f\left( x \right)$ is unknown function and
$$
_{C}D_{0t}^{\alpha }g=\left\{ \begin{aligned}
  & \frac{1}{\Gamma \left( 1-\alpha  \right)}\int\limits_{0}^{t}{\frac{{g}'\left( z \right)}{{{\left( t-z \right)}^{\alpha }}}dz,\,}\,\,0<\alpha <1, \\
 & \frac{dg}{dt},\,\,\,\,\,\,\,\,\,\,\,\,\,\,\,\,\,\,\,\,\,\,\,\,\,\,\,\,\,\,\,\,\,\,\,\,\,\,\,\,\,\,\,\,\alpha =1, \\
\end{aligned} \right.
$$
$$_{C}D_{t0}^{\beta }g=\left\{ \begin{aligned}
  & \frac{1}{\Gamma \left( 2-\beta  \right)}\int\limits_{t}^{0}{\frac{{g}''\left( z \right)}{{{\left( z-t \right)}^{\beta -1}}}dz,\,}\,\,1<\beta <2, \\
 & \frac{{{d}^{2}}g}{d{{t}^{2}}},\,\,\,\,\,\,\,\,\,\,\,\,\,\,\,\,\,\,\,\,\,\,\,\,\,\,\,\,\,\,\,\,\,\,\,\,\,\,\,\,\,\,\,\,\beta =2 \\
\end{aligned} \right.
$$
are fractional differential operators in a sense of Caputo (see (2.14)).

Here we are considering an inverse source problem for the equation (3.1), which is formulated as follows:

\textbf{Problem.} \emph{Find a pair of functions} $\left\{ u\left( x,t \right),\,f\left( t \right) \right\}$, \emph{satisfying}

i)	$u\left( x,t \right)\in C\left( \overline{\Omega } \right)\cap C_{x}^{2}\left( {{\Omega }^{+}}\cup {{\Omega }^{-}} \right),{{\,}_{C}}D_{0t}^{\alpha }u\in C\left( {{\Omega }^{+}} \right),{{\,}_{C}}D_{t0}^{\beta }u\in C\left( {{\Omega }^{-}} \right),$ \\ $f\left( t \right)\in C\left( 0,1 \right);$

ii)	\emph{equation (3.1) in} ${{\Omega }^{+}}$ and ${{\Omega }^{-}}$;

iii) \emph{the boundary conditions}
$$
u\left( 0,t \right)=u\left( 1,t \right),\,\,\,\,\,{{u}_{x}}\left( 0,t \right)=0,\,\,\,\,\,\,\,-p\le t\le q,\eqno (3.2)
$$
$$
u\left( x,-p \right)=\psi \left( x \right),\,\,\,\,\,u\left( x,q \right)=\varphi \left( x \right),\,\,\,\,0\le x\le 1;\eqno (3.3)
$$

iv)	\emph{the transmitting condition}
$$
\underset{t\to +0}{\mathop{\lim }}\,\,{{\,}_{C}}D_{0t}^{\alpha }u\left( x,t \right)=\underset{t\to -0}{\mathop{\lim }}\,\,{{\,}_{C}}D_{t0}^{\gamma }u\left( x,t \right),\,\,\,\,\,\,0<x<1,\eqno (3.4)
$$
where ${{\Omega }^{+}}=\Omega \cap \left\{ t>0 \right\},\,{{\Omega }^{-}}=\Omega \cap \left\{ t<0 \right\},\,\gamma =const\in \left( 0,1 \right],\,\psi \left( x \right),\,\varphi \left( x \right)$ are given functions such that $\varphi(0)=\varphi(1)$, $\psi(0)=\psi(1)$.

Using the method of separation of variables (Fourier's method) leads to the spectral problem (2.16) in the space variable $x$. Based on the properties of bi-orthogonal system, given in subsection 2.4, we look for a solution of problem (3.1)-(3.4) as follows
$$
u(x,t)={{V}_{0}}(t)+\sum\limits_{k=1}^{\infty }{}{{V}_{1k}}(t)\cos 2k\pi x+\sum\limits_{k=1}^{\infty }{}{{V}_{2k}}(t)\, x\sin 2k\pi x,\,\,\,t\ge 0,\eqno (3.5)
$$
$$
u(x,t)={{W}_{0}}(t)+\sum\limits_{k=1}^{\infty }{}{{W}_{1k}}(t)\cos 2k\pi x+\sum\limits_{k=1}^{\infty }{}{{W}_{2k}}(t)\, x\sin 2k\pi x,\,\,\,t\le 0,\eqno (3.6)
$$
$$
f(x)={{f}_{0}}+\sum\limits_{k=1}^{\infty }{}{{f}_{1k}}\,\cos 2k\pi x+\sum\limits_{k=1}^{\infty }{}{{f}_{2k}}\, x\sin 2k\pi x,\eqno (3.7)
$$
where the unknown coefficients are defined by
$$
\begin{aligned}
  & {{V}_{0}}(t)=2\int\limits_{0}^{1}{}u(x,t)(1-x)\,dx,\,\,\,\,\,\,t\ge 0,\\
 & {{V}_{1k}}(t)=4\int\limits_{0}^{1}{}u(x,t)(1-x)\cos 2k\pi x\,dx,\,\,\,\,\,\,t\ge 0, \\
 & {{V}_{2k}}(t)=4\int\limits_{0}^{1}{}u(x,t)\sin 2k\pi x\,dx,\,\,\,\,\,\,t\ge 0, \\
 & {{W}_{0}}(t)=2\int\limits_{0}^{1}{}u(x,t)(1-x)\,dx,\,\,\,\,\,\,\,t\le 0,\\
 & {{W}_{1k}}(t)=4\int\limits_{0}^{1}{}u(x,t)(1-x)\cos 2k\pi x\,dx,\,\,\,\,\,\,\,t\le 0, \\
 & {{W}_{2k}}(t)=4\int\limits_{0}^{1}{}u(x,t)\sin 2k\pi x\,dx,\,\,\,\,\,\,\,t\le 0, \\
 & {{f}_{0}}=2\int\limits_{0}^{1}{}f(x)(1-x)\,dx,\\
 & {{f}_{1k}}=4\int\limits_{0}^{1}{}f(x)(1-x)\cos 2k\pi x\,dx,\\
 &{{f}_{2k}}=4\int\limits_{0}^{1}{}f(x)\sin 2k\pi x\,dx. \\
\end{aligned}\eqno (3.8)
$$

Based on (3.8), we introduce the following functions, in which $\varepsilon >0$ is sufficiently small number:

\newpage

$$
\begin{aligned}
  & {{V}_{0,\varepsilon }}(t)=2\int\limits_{\varepsilon }^{1-\varepsilon }{}u(x,t)(1-x)\,dx,\,\,\,\,t\ge 0,\\
  & {{V}_{1k,\varepsilon }}(t)=4\int\limits_{\varepsilon }^{1-\varepsilon }{}u(x,t)(1-x)\cos 2k\pi x\,dx, \,\,\,\,t\ge 0,\\
 & {{V}_{2k,\varepsilon }}(t)=4\int\limits_{\varepsilon }^{1-\varepsilon }{}u(x,t)\sin 2k\pi x\,dx,\,\,\,\,t\ge 0, \\
 & {{W}_{0,\varepsilon }}(t)=2\int\limits_{\varepsilon }^{1-\varepsilon }{}u(x,t)(1-x)\,dx,\,\,\,t\le 0,\\
 & {{W}_{1k,\varepsilon }}(t)=4\int\limits_{\varepsilon }^{1-\varepsilon }{}u(x,t)(1-x)\cos 2k\pi x\,dx,\,\,\,t\le 0, \\
 & {{W}_{2k,\varepsilon }}(t)=4\int\limits_{\varepsilon }^{1-\varepsilon }{}u(x,t)\sin 2k\pi x\,dx,\,\,\,t\le 0, \\
 & {{f}_{0,\varepsilon }}=2\int\limits_{\varepsilon }^{1-\varepsilon }{}f(x)(1-x)d\,x,\\
 & {{f}_{1k,\varepsilon }}=4\int\limits_{\varepsilon }^{1-\varepsilon }{}f(x)(1-x)\cos 2k\pi x\,dx,\\
 &{{f}_{2k,\varepsilon }}=4\int\limits_{\varepsilon }^{1-\varepsilon }{}f(x)\sin 2k\pi x\,dx. \\
\end{aligned}\eqno (3.9)
$$

 Applying the operators $_{C}D_{0t}^{\alpha }\left( \cdot  \right)$ for $t\in \left( 0,q \right)$ and $_{C}D_{t0}^{\beta }\left( \cdot  \right)$ for $t\in \left( -p,0 \right)$ , and using (3.1), we get

\newpage

$$
\begin{aligned}
  & _{C}D_{0t}^{\alpha }{{V}_{0,\varepsilon }}\left( t \right)=2\int\limits_{\varepsilon }^{1-\varepsilon }{{{u}_{xx}}\left( x,t \right)\left( 1-x \right)\,dx}+{{f}_{0,\varepsilon }},\,\,\,\,\,t>0, \\
 & _{C}D_{0t}^{\alpha }{{V}_{1k,\varepsilon }}\left( t \right)=4\int\limits_{\varepsilon }^{1-\varepsilon }{{{u}_{xx}}\left( x,t \right)\left( 1-x \right)\cos 2k\pi x\,dx}+{{f}_{1k,\varepsilon }},\,\,\,\,\,t>0, \\
 & _{C}D_{0t}^{\alpha }{{V}_{2k,\varepsilon }}\left( t \right)=4\int\limits_{\varepsilon }^{1-\varepsilon }{{{u}_{xx}}\left( x,t \right)\sin 2k\pi x\,dx}+{{f}_{2k,\varepsilon }},\,\,\,\,\,t>0, \\
 & _{C}D_{t0}^{\beta }{{W}_{0,\varepsilon }}\left( t \right)=2\int\limits_{\varepsilon }^{1-\varepsilon }{{{u}_{xx}}\left( x,t \right)\left( 1-x \right)\,dx}+{{f}_{0,\varepsilon }},\,\,\,\,\,t<0, \\
 & _{C}D_{t0}^{\beta }{{W}_{1k,\varepsilon }}\left( t \right)=4\int\limits_{\varepsilon }^{1-\varepsilon }{{{u}_{xx}}\left( x,t \right)\left( 1-x \right)\cos 2k\pi x\,dx}+{{f}_{1k,\varepsilon }},\,\,\,\,\,t<0, \\
 & _{C}D_{t0}^{\beta }{{W}_{2k,\varepsilon }}\left( t \right)=4\int\limits_{\varepsilon }^{1-\varepsilon }{{{u}_{xx}}\left( x,t \right)\sin 2k\pi x\,dx}+{{f}_{2k,\varepsilon }},\,\,\,\,\,t<0. \\
\end{aligned}\eqno (3.10)
$$
On integrating by parts and taking the limit as $\varepsilon \to 0$, we obtain the following set of equations for the unknown coefficients
$$
_{C}D_{0t}^{\alpha }{{V}_{0}}\left( t \right)={{f}_{0}},\,\,t\ge 0,\eqno (3.11)
$$
$$
_{C}D_{t0}^{\beta }{{W}_{0}}\left( t \right)={{f}_{0}},\,\,t<0,\eqno (3.12)
$$
$$
_{C}D_{0t}^{\alpha }{{V}_{1k}}\left( t \right)+{{\left( 2k\pi  \right)}^{2}}{{V}_{1k}}\left( t \right)={{f}_{1k}}+4k\pi {{V}_{2k}}\left( t \right),\,\,t\ge 0,\eqno (3.13)
$$
$$
_{C}D_{t0}^{\beta }{{W}_{1k}}\left( t \right)+{{\left( 2k\pi  \right)}^{2}}{{W}_{1k}}\left( t \right)={{f}_{1k}}+4k\pi {{W}_{2k}}\left( t \right),\,\,t<0,\eqno (3.14)
$$
$$
_{C}D_{0t}^{\alpha }{{V}_{2k}}\left( t \right)+{{\left( 2k\pi  \right)}^{2}}{{V}_{2k}}\left( t \right)={{f}_{2k}}\,\,,\,\,t\ge 0,\eqno (3.15)
$$
$$
_{C}D_{t0}^{\beta }{{W}_{2k}}\left( t \right)+{{\left( 2k\pi  \right)}^{2}}{{W}_{2k}}\left( t \right)={{f}_{2k}},\,\,t<0.\eqno (3.16)
$$
General solutions of (3.11) and (3.15) can be written as follows (see [17, p.17])
$$
{{V}_{0}}(t)={{V}_{0}}(0)+\frac{{{f}_{0}}}{\Gamma \left( \alpha +1 \right)}{{t}^{\alpha }},\eqno (3.17)
$$
$$
{{V}_{2k}}(t)={{V}_{2k}}(0){{E}_{\alpha ,1}}\left( -{{(2k\pi )}^{2}}{{t}^{\alpha }} \right)+{{f}_{2k}}{{t}^{\alpha }}{{E}_{\alpha ,\alpha +1}}\left( -{{(2k\pi )}^{2}}{{t}^{\alpha }} \right),\eqno (3.18)
$$
where ${{E}_{\alpha ,\beta }}\left( z \right)$ is the two parameter Mittag-Leffler function, defined by (2.4).

Hence, a general solution of (3.13) is given by (see Appendix A1 for details)
$$
\begin{aligned}
  & {{V}_{1k}}(t)={{V}_{1k}}(0){{E}_{\alpha ,1}}\left( -{{(2k\pi )}^{2}}{{t}^{\alpha }} \right)+{{f}_{1k}}\cdot {{t}^{\alpha }}\cdot {{E}_{\alpha ,\alpha +1}}\left( -{{(2k\pi )}^{2}}{{t}^{\alpha }} \right)+ \\
 & +4k\pi \cdot {{V}_{2k}}(0)\cdot {{t}^{\alpha }}\cdot {{E}_{1}}\left( \begin{matrix}
   1,1;1,1 & |-{{(2k\pi )}^{2}}{{t}^{\alpha }}  \\
   \alpha +1,\alpha ,\alpha ;1,1;1,1 & |-{{(2k\pi )}^{2}}{{t}^{\alpha }}  \\
\end{matrix} \right)+ \\
 & +4k\pi \cdot {{f}_{2k}}\cdot {{t}^{2\alpha }}\cdot {{E}_{1}}\left( \begin{matrix}
   1,1;1,1 & |-{{(2k\pi )}^{2}}{{t}^{\alpha }}  \\
   2\alpha +1,\alpha ,\alpha ;1,1;1,1 & |-{{(2k\pi )}^{2}}{{t}^{\alpha }}  \\
\end{matrix} \right), \\
\end{aligned}\eqno (3.19)
$$
where
$${{E}_{1}}\left( \begin{matrix}
   {{\gamma }_{1}},{{\alpha }_{1}};{{\gamma }_{2}},{{\beta }_{1}} & |\,\,x  \\
   {{\delta }_{1}},{{\alpha }_{2}},{{\beta }_{2}};{{\delta }_{2}},{{\alpha }_{3}};{{\delta }_{3}},{{\beta }_{3}} & |\,\,y  \\
\end{matrix} \right)
$$ is the Mittag-Leffler type function of two variables, defined by (2.7).

Similarly, general solutions of (3.12), (3.14) and (3.16) are given by (See Appendix A2 for integral forms of solutions as in [17, p.17])
$$
{{W}_{0}}(t)={{W}_{0}}(0)-t{{W}_{0}}^{\prime }(0)+\frac{{{f}_{0}}}{\Gamma \left( \beta +1 \right)}{{(-t)}^{\beta }},\eqno (3.20)
$$
$$
\begin{aligned}
  & {{W}_{1k}}(t)={{W}_{1k}}(0){{E}_{\beta ,1}}\left( -{{(2k\pi )}^{2}}{{(-t)}^{\beta }} \right)-t{{W}_{1k}}^{\prime }(0){{E}_{\beta ,2}}\left( -{{(2k\pi )}^{2}}{{(-t)}^{\beta }} \right)+ \\
 & +{{f}_{1k}}\cdot {{(-t)}^{\beta }}{{E}_{\beta ,\beta +1}}\left( -{{(2k\pi )}^{2}}{{(-t)}^{\beta }} \right)+ \\
 & +4k\pi \cdot {{W}_{2k}}(0)\cdot {{(-t)}^{\beta }}{{E}_{1}}\left( \begin{matrix}
   1,1;1,1 & |-{{(2k\pi )}^{2}}{{\left( -t \right)}^{\beta }}  \\
   \beta +1,\beta ,\beta ;1,1;1,1 & |-{{(2k\pi )}^{2}}{{\left( -t \right)}^{\beta }}  \\
\end{matrix} \right)+ \\
 & +4k\pi \cdot {{W}_{2k}}^{\prime }(0)\cdot {{(-t)}^{\beta +1}}{{E}_{1}}\left( \begin{matrix}
   1,1;1,1 & |-{{(2k\pi )}^{2}}{{\left( -t \right)}^{\beta }}  \\
   \beta +2,\beta ,\beta ;1,1;1,1 & |-{{(2k\pi )}^{2}}{{\left( -t \right)}^{\beta }}  \\
\end{matrix} \right)+ \\
 & +4k\pi \cdot {{f}_{2k}}\cdot {{(-t)}^{2\beta }}{{E}_{1}}\left( \begin{matrix}
   1,1;1,1 & |-{{(2k\pi )}^{2}}{{\left( -t \right)}^{\beta }}  \\
   2\beta +1,\beta ,\beta ;1,1;1,1 & |-{{(2k\pi )}^{2}}{{\left( -t \right)}^{\beta }}  \\
\end{matrix} \right), \\
\end{aligned}\eqno (3.21)
$$
$$
\begin{aligned}
  & {{W}_{2k}}(t)={{W}_{2k}}(0){{E}_{\beta ,1}}\left( -{{(2k\pi )}^{2}}{{(-t)}^{\beta }} \right)-t{{W}_{2k}}^{\prime }(0){{E}_{\beta ,2}}\left( -{{(2k\pi )}^{2}}{{(-t)}^{\beta }} \right)+ \\
 & +{{f}_{2k}}{{(-t)}^{\beta }}{{E}_{\beta ,\beta +1}}\left( -{{(2k\pi )}^{2}}{{(-t)}^{\beta }} \right), \\
\end{aligned}\eqno (3.22)
$$
where the unknown constants ${{f}_{0}},\,{{f}_{1k,}}{{f}_{2k,}}\,{{V}_{0}}\left( 0 \right),\,{{V}_{1k}}\left( 0 \right),\,{{V}_{2k}}\left( 0 \right),\,{{W}_{0}}\left( 0 \right),\,{{W}_{1k}}\left( 0 \right),$ ${{W}_{2k}}\left( 0 \right),\,{{W}_{1k}}^{\prime }\left( 0 \right),$ ${{W}_{2k}}^{\prime }\left( 0 \right)$  to be determined using the given boundary and transmitting conditions.

Starting with conditions (3.3), we get
$$
{{W}_{0}}(0)+p{{W}_{0}}^{\prime }(0)+\frac{{{f}_{0}}}{\Gamma \left( \beta +1 \right)}{{p}^{\beta }}={{\psi }_{0}},\eqno (3.23)
$$
$$
\begin{aligned}
  & {{W}_{1k}}(0){{E}_{\beta ,1}}\left( -{{(2k\pi )}^{2}}{{p}^{\beta }} \right)+p{{W}_{1k}}^{\prime }(0){{E}_{\beta ,2}}\left( -{{(2k\pi )}^{2}}{{p}^{\beta }} \right)+{{f}_{1k}}{{p}^{\beta }}{{E}_{\beta ,\beta +1}}\left( -{{(2k\pi )}^{2}}{{p}^{\beta }} \right)+ \\
 & +4k\pi \cdot {{W}_{2k}}(0)\cdot {{p}^{\beta }}{{E}_{1}}\left( \begin{matrix}
   1,1;1,1 & |-{{(2k\pi )}^{2}}{{p}^{\beta }}  \\
   \beta +1,\beta ,\beta ;1,1;1,1 & |-{{(2k\pi )}^{2}}{{p}^{\beta }}  \\
\end{matrix} \right)+ \\
 & +4k\pi \cdot {{W}_{2k}}^{\prime }(0)\cdot {{p}^{\beta +1}}{{E}_{1}}\left( \begin{matrix}
   1,1;1,1 & |-{{(2k\pi )}^{2}}{{p}^{\beta }}  \\
   \beta +2,\beta ,\beta ;1,1;1,1 & |-{{(2k\pi )}^{2}}{{p}^{\beta }}  \\
\end{matrix} \right)+ \\
 & +4k\pi \cdot {{f}_{2k}}\cdot {{p}^{2\beta }}{{E}_{1}}\left( \begin{matrix}
   1,1;1,1 & |-{{(2k\pi )}^{2}}{{p}^{\beta }}  \\
   2\beta +1,\beta ,\beta ;1,1;1,1 & |-{{(2k\pi )}^{2}}{{p}^{\beta }}  \\
\end{matrix} \right)={{\psi }_{1k}}, \\
\end{aligned}\eqno (3.24)
$$
$$
\begin{aligned}
& {{W}_{2k}}(0){{E}_{\beta ,1}}\left( -{{(2k\pi )}^{2}}{{p}^{\beta }} \right)+p{{W}_{2k}}'(0){{E}_{\beta ,2}}\left( -{{(2k\pi )}^{2}}{{p}^{\beta }} \right)+\\
& +{{f}_{2k}}{{p}^{\beta }}{{E}_{\beta ,\beta +1}}\left( -{{(2k\pi )}^{2}}{{p}^{\beta }} \right)={{\psi }_{2k}},\\
\end{aligned}\eqno (3.25)
$$
$$
{{V}_{0}}(0)+\frac{{{f}_{0}}}{\Gamma \left( \alpha +1 \right)}{{q}^{\alpha }}={{\varphi }_{0}},\eqno (3.26)
$$
$$
\begin{aligned}
  & {{V}_{1k}}(0){{E}_{\alpha ,1}}\left( -{{(2k\pi )}^{2}}{{q}^{\alpha }} \right)+{{f}_{1k}}{{q}^{\alpha }}{{E}_{\alpha ,\alpha +1}}\left( -{{(2k\pi )}^{2}}{{q}^{\alpha }} \right)+ \\
 & +4k\pi \cdot {{V}_{2k}}(0)\cdot {{q}^{\alpha }}{{E}_{1}}\left( \begin{matrix}
   1,1;1,1 & |-{{(2k\pi )}^{2}}{{q}^{\alpha }}  \\
   \alpha +1,\alpha ,\alpha ;1,1;1,1 & |-{{(2k\pi )}^{2}}{{q}^{\alpha }}  \\
\end{matrix} \right)+ \\
 & +4k\pi \cdot {{f}_{2k}}\cdot {{q}^{2\alpha }}{{E}_{1}}\left( \begin{matrix}
   1,1;1,1 & |-{{(2k\pi )}^{2}}{{q}^{\alpha }}  \\
   2\alpha +1,\alpha ,\alpha ;1,1;1,1 & |-{{(2k\pi )}^{2}}{{q}^{\alpha }}  \\
\end{matrix} \right)={{\varphi }_{1k}}, \\
\end{aligned}\eqno (3.27)
$$
$$
{{V}_{2k}}(0){{E}_{\alpha ,1}}\left( -{{(2k\pi )}^{2}}{{q}^{\alpha }} \right)+{{f}_{2k}}{{q}^{\alpha }}{{E}_{\alpha ,\alpha +1}}\left( -{{(2k\pi )}^{2}}{{q}^{\alpha }} \right)={{\varphi }_{2k}}, \eqno (3.28)
$$
where
$$
{{\varphi }_{0}}=2\int\limits_{0}^{1}{}\varphi (x)(1-x)\,dx,{{\varphi }_{1k}}=4\int\limits_{0}^{1}{}\varphi (x)(1-x)\cos 2k\pi x\,dx,{{\varphi }_{2k}}=4\int\limits_{0}^{1}{}\varphi (x)\sin 2k\pi x\,dx,
$$
$$
{{\psi }_{0}}=2\int\limits_{0}^{1}{}\psi (x)(1-x)\,dx,{{\psi }_{1k}}=4\int\limits_{0}^{1}{}\psi (x)(1-x)\cos 2k\pi x\,dx,{{\psi }_{2k}}=4\int\limits_{0}^{1}{}\psi (x)\sin 2k\pi x\,dx.
$$

Next using $u\left( x,+0 \right)=u\left( x,-0 \right)$ (due to $u\left( x,t \right)\in C\left( \overline{\Omega } \right)$), we get
$$
{{V}_{0}}(0)={{W}_{0}}(0),\,\,\,\,\,\,{{V}_{1k}}(0)={{W}_{1k}}(0),\,\,\,\,\,\,{{V}_{2k}}(0)={{W}_{2k}}(0).\eqno (3.29)
$$

Now, to use the transmitting condition (3.4), we first take the limit as $t\rightarrow +0$ in equations (3.11), (3.13) and (3.15) to get
$$
\begin{aligned}
  & \lim\limits_{t\to +0}\,{{\,}_{C}}D_{0t}^{\alpha }{{V}_{0}}(t)={{f}_{0}},\,\,\lim\limits_{t\to +0}\,{{\,}_{C}}D_{0t}^{\alpha }{{V}_{1k}}(t)={{f}_{1k}}+4k\pi {{V}_{2k}}(0)-{{\left( 2k\pi  \right)}^{2}}{{V}_{1k}}\left( 0 \right),\,\,\, \\
 & \lim\limits_{t\to +0}\,{{\,}_{C}}D_{0t}^{\alpha }{{V}_{2k}}(t)={{f}_{2k}}-{{(2k\pi )}^{2}}{{V}_{2k}}(0), \\
\end{aligned}\eqno (3.30)
$$
and from (3.20) we obtain
$$
_{C}D_{t0}^{\gamma }{{W}_{0}}\left( t \right)={{\,}_{RL}}D_{t0}^{\gamma }{{W}_{0}}\left( t \right)-\frac{{{W}_{0}}\left( 0 \right)}{\Gamma \left( 1-\gamma  \right){{\left| t \right|}^{\gamma }}},
$$
where
$$
_{RL}D_{t0}^{\gamma }{{W}_{0}}\left( t \right)={{\,}_{RL}}D_{t0}^{\gamma }{{W}_{0}}\left( 0 \right)+{{\,}_{RL}}D_{t0}^{\gamma }\left[ -t{{W}_{0}}^{\prime }\left( 0 \right) \right]+{{\,}_{RL}}D_{t0}^{\gamma }\left[ \frac{{{f}_{0}}}{\Gamma \left( \beta +1 \right)}{{\left( -t \right)}^{\beta }} \right].
$$
The three terms on the right hand side of the above equation can be determined using the definition of Riemann-Liouville fractional differential operator (2.13), Euler's beta-function (2.3) and properties of Euler's gamma-function (2.1), (2.2). For example, the last term can be written as
$$
\begin{aligned}
  & _{RL}D_{t0}^{\gamma }\left[ \frac{{{f}_{0}}}{\Gamma \left( \beta +1 \right)}{{\left( -t \right)}^{\beta }} \right]=-\frac{{{f}_{0}}}{\Gamma \left( \beta +1 \right)\,\Gamma \left( 1-\gamma  \right)}\cdot \frac{d}{dt}\int\limits_{t}^{0}{\frac{{{\left( -\xi  \right)}^{\beta }}}{{{\left( \xi -t \right)}^{\gamma }}}d\xi }=\\
  &=\frac{{{f}_{0}}}{\Gamma \left( \beta +1 \right)\,\Gamma \left( 1-\gamma  \right)}\cdot \frac{d}{d\left( -t \right)}\left( {{\left( -t \right)}^{\beta +1-\gamma }} \right)
 \times \frac{\Gamma \left( \beta +1 \right)\Gamma \left( 1-\gamma  \right)}{\Gamma \left( \beta +2-\gamma  \right)}=\\
 &=\frac{{{f}_{0}}\,\Gamma \left( \beta +1 \right)}{\,\Gamma \left( \beta +1-\gamma  \right)}{{\left( -t \right)}^{\beta -\gamma }}. \\
\end{aligned}
$$
Hence, we deduce
$$
_{C}D_{t0}^{\gamma }{{W}_{0}}\left( t \right)=\frac{{{W}_{0}}^{\prime }\left( 0 \right)}{\Gamma \left( 2-\gamma  \right)}{{\left( -t \right)}^{1-\gamma }}+\frac{{{f}_{0}}}{\,\Gamma \left( \beta +1-\gamma  \right)\,}{{\left( -t \right)}^{\beta -\gamma }}.\eqno (3.31)
$$
The expression for $_{C}D_{t0}^{\gamma }{{W}_{2k}}\left( t \right)$ can be obtained using (2.15) and (3.22):
$$
\begin{aligned}
  & _{C}D_{t0}^{\gamma }{{W}_{2k}}\left( t \right)={{\,}_{RL}}D_{t0}^{\gamma }\left[ {{W}_{2k}}\left( 0 \right){{E}_{\beta ,1}}\left( -{{\left( 2k\pi  \right)}^{2}}{{\left( -t \right)}^{\beta }} \right) \right]-\\
  &-{{\,}_{RL}}D_{t0}^{\gamma }\left[ t{{W}_{2k}}^{\prime }\left( 0 \right){{E}_{\beta ,2}}\left( -{{\left( 2k\pi  \right)}^{2}}{{\left( -t \right)}^{\beta }} \right) \right]+ \\
 & +{{\,}_{RL}}D_{t0}^{\gamma }\left[ {{f}_{2k}}{{\left( -t \right)}^{\beta }}{{E}_{\beta ,\beta +1}}\left( -{{\left( 2k\pi  \right)}^{2}}{{\left( -t \right)}^{\beta }} \right) \right]-\frac{{{W}_{2k}}\left( 0 \right)}{\Gamma \left( 1-\gamma  \right)}{{\left( -t \right)}^{-\gamma }}. \\
\end{aligned}
$$
The three Riemann-Liouville derivative terms on the right hand side of the above equation can be evaluated using the differentiation formula (2.6). For example, the first term can be evaluated by setting $\alpha =\beta ,\,\overline{\beta }=1,\,k=0,\,\lambda =-{{\left( 2k\pi  \right)}^{2}},\,y=-t$:
$$
_{RL}D_{t0}^{\gamma }\left[ {{W}_{2k}}\left( 0 \right){{E}_{\beta ,1}}\left( -{{\left( 2k\pi  \right)}^{2}}{{\left( -t \right)}^{\beta }} \right) \right]={{W}_{2k}}\left( 0 \right){{\left( -t \right)}^{-\gamma }}{{E}_{\beta ,1-\gamma }}\left( -{{\left( 2k\pi  \right)}^{2}}{{\left( -t \right)}^{\beta }} \right).
$$
Similarly, one can evaluate the other two terms. Substituting back and using the property (2.5), the expression of $ _{C}D_{t0}^{\gamma }{{W}_{2k}}\left( t \right)$ can be rewritten as
$$
\begin{aligned}
  & _{C}D_{t0}^{\gamma }{{W}_{2k}}\left( t \right)=\left[ {{f}_{2k}}-{{W}_{2k}}\left( 0 \right) \right]{{\left( 2k\pi  \right)}^{2}}{{\left( -t \right)}^{\beta -\gamma }}{{E}_{\beta ,\beta +1-\gamma }}\left( -{{\left( 2k\pi  \right)}^{2}}{{\left( -t \right)}^{\beta }} \right)+ \\
 & +{{W}_{2k}}^{\prime }\left( 0 \right){{\left( -t \right)}^{1-\gamma }}{{E}_{\beta ,2-\gamma }}\left( -{{\left( 2k\pi  \right)}^{2}}{{\left( -t \right)}^{\beta }} \right). \\
\end{aligned}\eqno (3.32)
$$
Now, using equation (3.21) and relation (2.15), we obtain the following expression for $_{C}D_{t0}^{\gamma }{{W}_{1k}}\left( t \right)$:
$$
\begin{aligned}
  & _{C}D_{t0}^{\gamma }{{W}_{1k}}\left( t \right)={{W}_{1k}}\left( 0 \right){{\,}_{RL}}D_{t0}^{\gamma }\left[ {{E}_{\beta ,1}}\left( -{{\left( 2k\pi  \right)}^{2}}{{\left( -t \right)}^{\beta }} \right) \right]+\\
  &+{{W}_{1k}}^{\prime }\left( 0 \right){{\,}_{RL}}D_{t0}^{\gamma }\left[ -t\cdot {{E}_{\beta ,2}}\left( -{{\left( 2k\pi  \right)}^{2}}{{\left( -t \right)}^{\beta }} \right) \right]+ \\
 & +{{f}_{1k}}\cdot {{\,}_{RL}}D_{t0}^{\gamma }\left[ {{\left( -t \right)}^{\beta }}{{E}_{\beta ,\beta +1}}\left( -{{\left( 2k\pi  \right)}^{2}}{{\left( -t \right)}^{\beta }} \right) \right]+ \\
 & +4k\pi {{W}_{2k}}\left( 0 \right)\cdot {{\,}_{RL}}D_{t0}^{\gamma }\left[ {{\left( -t \right)}^{\beta }}{{E}_{1}}\left( \begin{matrix}
   1,1;1,1 & |\,\,-{{\left( 2k\pi  \right)}^{2}}{{\left( -t \right)}^{\beta }}  \\
   \beta +1,\beta ,\beta ;1,1;1,1 & |\,\,-{{\left( 2k\pi  \right)}^{2}}{{\left( -t \right)}^{\beta }}  \\
\end{matrix} \right) \right]+ \\
 & +4k\pi \cdot {{W}_{2k}}^{\prime }\left( 0 \right)\cdot {{\,}_{RL}}D_{t0}^{\gamma }\left[ {{\left( -t \right)}^{\beta +1}}{{E}_{1}}\left( \begin{matrix}
   1,1;1,1 & |\,\,-{{\left( 2k\pi  \right)}^{2}}{{\left( -t \right)}^{\beta }}  \\
   \beta +2,\beta ,\beta ;1,1;1,1 & |\,\,-{{\left( 2k\pi  \right)}^{2}}{{\left( -t \right)}^{\beta }}  \\
\end{matrix} \right) \right]+ \\
 & +4k\pi \cdot {{f}_{2k}}\cdot {{\,}_{RL}}D_{t0}^{\gamma }\left[ {{\left( -t \right)}^{2\beta }}{{E}_{1}}\left( \begin{matrix}
   1,1;1,1 & |\,\,-{{\left( 2k\pi  \right)}^{2}}{{\left( -t \right)}^{\beta }}  \\
   2\beta +1,\beta ,\beta ;1,1;1,1 & |\,\,-{{\left( 2k\pi  \right)}^{2}}{{\left( -t \right)}^{\beta }}  \\
\end{matrix} \right) \right]-\\
&-\frac{{{W}_{1k}}\left( 0 \right)}{\Gamma \left( 1-\gamma  \right)}{{\left( -t \right)}^{-\gamma }}. \\
\end{aligned}
$$
This expression can be simplified by calculating the Riemann-Liouville derivative terms on the right hand side of the above equation using the property (2.5) and formulas (2.6), (2.9). For example, the last Riemann-Liouville derivative term can be calculated by setting
${{\delta }_{1}}=2\beta +1,\,{{\gamma }_{1}}={{\gamma }_{2}}={{\alpha }_{1}}={{\beta }_{1}}={{\delta }_{2}}={{\delta }_{3}}={{\alpha }_{3}}={{\beta }_{3}}=1,\,{{\alpha }_{2}}={{\beta }_{2}}=\beta ,y=t,\,{{\omega }_{1}}={{\omega }_{2}}=-{{\left( 2k\pi  \right)}^{2}}$
in (2.9) to obtain
$$
\begin{aligned}
  & _{RL}D_{t0}^{\gamma }\left[ {{\left( -t \right)}^{2\beta }}{{E}_{1}}\left( \begin{matrix}
   1,1;1,1 & |\,\,-{{\left( 2k\pi  \right)}^{2}}{{\left( -t \right)}^{\beta }}  \\
   2\beta +1,\beta ,\beta ;1,1;1,1 & |\,\,-{{\left( 2k\pi  \right)}^{2}}{{\left( -t \right)}^{\beta }}  \\
\end{matrix} \right) \right]= \\
 & ={{\left( -t \right)}^{2\beta -\gamma }}{{E}_{1}}\left( \begin{matrix}
   1,1;1,1 & |\,\,-{{\left( 2k\pi  \right)}^{2}}{{\left( -t \right)}^{\beta }}  \\
   2\beta -\gamma +1,\beta ,\beta ;1,1;1,1 & |\,\,-{{\left( 2k\pi  \right)}^{2}}{{\left( -t \right)}^{\beta }}  \\
\end{matrix} \right). \\
\end{aligned}
$$
Hence, we get
$$
\begin{aligned}
  & _{C}D_{t0}^{\gamma }{{W}_{1k}}\left( t \right)=\left[ {{f}_{1k}}-{{\left( 2k\pi  \right)}^{2}}{{W}_{1k}}\left( 0 \right) \right]{{\left( -t \right)}^{\beta -\gamma }}{{E}_{\beta ,1+\beta -\gamma }}\left( -{{\left( 2k\pi  \right)}^{2}}{{\left( -t \right)}^{\beta }} \right)+ \\
 & +{{W}_{1k}}^{\prime }\left( 0 \right){{\left( -t \right)}^{1-\gamma }}{{E}_{\beta ,2-\gamma }}\left( -{{\left( 2k\pi  \right)}^{2}}{{\left( -t \right)}^{\beta }} \right)+ \\
 & +4k\pi {{W}_{2k}}\left( 0 \right){{\left( -t \right)}^{\beta -\gamma }}{{E}_{1}}\left( \begin{matrix}
   1,1;1,1 & |\,\,-{{\left( 2k\pi  \right)}^{2}}{{\left( -t \right)}^{\beta }}  \\
   \beta +1-\gamma ,\beta ,\beta ;1,1;1,1 & |\,\,-{{\left( 2k\pi  \right)}^{2}}{{\left( -t \right)}^{\beta }}  \\
\end{matrix} \right)+ \\
 & +4k\pi {{W}_{2k}}^{\prime }\left( 0 \right){{\left( -t \right)}^{1+\beta -\gamma }}{{E}_{1}}\left( \begin{matrix}
   1,1;1,1 & |\,\,-{{\left( 2k\pi  \right)}^{2}}{{\left( -t \right)}^{\beta }}  \\
   2+\beta -\gamma ,\beta ,\beta ;1,1;1,1 & |\,\,-{{\left( 2k\pi  \right)}^{2}}{{\left( -t \right)}^{\beta }}  \\
\end{matrix} \right)+ \\
 & +4k\pi \cdot {{f}_{2k}}\cdot {{\left( -t \right)}^{2\beta -\gamma }}{{E}_{1}}\left( \begin{matrix}
   1,1;1,1 & |\,\,-{{\left( 2k\pi  \right)}^{2}}{{\left( -t \right)}^{\beta }}  \\
   2\beta +1-\gamma ,\beta ,\beta ;1,1;1,1 & |\,\,-{{\left( 2k\pi  \right)}^{2}}{{\left( -t \right)}^{\beta }}  \\
\end{matrix} \right). \\
\end{aligned}\eqno (3.33)
$$

Finally, applying the transmitting condition (3.4), taking into account (3.30)-(3.33), we obtain
$$
{{f}_{0}}=0,\,{{f}_{1k}}+4k\pi {{V}_{2k}}\left( 0 \right)-{{\left( 2k\pi  \right)}^{2}}{{V}_{1k}}\left( 0 \right)=0,\,\,{{f}_{2k}}-{{\left( 2k\pi  \right)}^{2}}{{V}_{2k}}\left( 0 \right)=0.\eqno (3.34)
$$
The solutions of algebraic equations (3.23)-(3.29) and (3.34) are given by
$$
\begin{aligned}
  & {{f}_{0}}=0,\\
 & {{f}_{1k}}={{\left( 2k\pi  \right)}^{2}}{{\varphi }_{1k}}-4k\pi {{\varphi }_{2k}},\\
 & {{f}_{2k}}={{\left( 2k\pi  \right)}^{2}}{{\varphi }_{2k}},\\
 &{{V}_{0}}\left( 0 \right)={{\varphi }_{0}},\\
 & {{V}_{1k}}\left( 0 \right)={{\varphi }_{1k}},\\
& {{V}_{2k}}\left( 0 \right)={{\varphi }_{2k}},\\
& {{W}_{0}}^{\prime }\left( 0 \right)=\frac{{{\psi }_{0}}-{{\varphi }_{0}}}{p},\\
& {{W}_{0}}\left( 0 \right)={{\varphi }_{0}},\\
& {{W}_{1k}}\left( 0 \right)={{\varphi }_{1k}}, \\
 & {{W}_{1k}}^{\prime }\left( 0 \right)={{\psi }_{1k}}-{{\varphi }_{1k}}+\frac{{{\psi }_{2k}}-{{\varphi }_{2k}}}{{{E}_{\beta ,2}}\left( -{{\left( 2k\pi  \right)}^{2}}{{p}^{\beta }} \right)}\times\\
 &\times 4k\pi {{p}^{\beta }}{{E}_{1}}\left( \begin{matrix}
   1,1;1,1 & |-{{(2k\pi )}^{2}}{{p}^{\beta }}  \\
   \beta +2,\beta ,\beta ;1,1;1,1 & |-{{(2k\pi )}^{2}}{{p}^{\beta }}  \\
\end{matrix} \right), \\
& {{W}_{2k}}\left( 0 \right)={{\varphi }_{2k}},\\
& {{W}_{2k}}^{\prime }\left( 0 \right)=\frac{{{\psi }_{2k}}-{{\varphi }_{2k}}}{p\,{{E}_{\beta ,2}}\left( -{{\left( 2k\pi  \right)}^{2}}{{p}^{\beta }} \right)}. \\
\end{aligned}\eqno (3.35)
$$

\section{The existence and uniqueness of the solution}

\subsection{The existence}

Substituting the obtained expressions in (3.35) into (3.17)-(3.22), the unknown coefficients are represented in terms of the given data and hence, the solution (3.5)-(3.7) can be rewritten as
$$
\begin{aligned}
&u(x,t)={{\varphi }_{0}}-\frac{{{\psi }_{0}}-{{\varphi }_{0}}}{p}\cdot t+\sum\limits_{k=1}^{\infty }{}{{W}_{1k}}(t)\cos 2k\pi x+\sum\limits_{k=1}^{\infty }{}{{W}_{2k}}(t)\, x\sin 2k\pi x,\,\,t\le 0,\\
&u(x,t)={{{\varphi }_{0}}}+\sum\limits_{k=1}^{\infty }{}{{\varphi }_{1k}} \cos 2k\pi x+\sum\limits_{k=1}^{\infty }{}{{\varphi }_{2k}}\, x\sin 2k\pi x,\,\,t\ge 0,\\
\end{aligned}\eqno (4.1)
$$
$$
f(x)=\sum\limits_{k=1}^{\infty }{}\left({{\left( 2k\pi  \right)}^{2}}{{\varphi }_{1k}}-4k\pi {{\varphi }_{2k}}\right)\,\cos 2k\pi x+\sum\limits_{k=1}^{\infty }{}{{\left( 2k\pi  \right)}^{2}}{{\varphi }_{2k}}\, x\sin 2k\pi x,\eqno (4.2)
$$
where
$$
\begin{aligned}
   & {{W}_{1k}}\left( t \right)={{\varphi }_{1k}}-\overline{{{\varphi }_{1k}}}\cdot t{{E}_{\beta ,2}}\left( -{{\left( 2k\pi  \right)}^{2}}{{\left( -t \right)}^{\beta }} \right)+\frac{{{\psi }_{2k}}-{{\varphi }_{2k}}}{p{{E}_{\beta ,2}}\left( -{{\left( 2k\pi  \right)}^{2}}{{p}^{\beta }} \right)}\times  \\
 & \times 4k\pi {{\left( -t \right)}^{\beta +1}}{{E}_{1}}\left( \begin{matrix}
   1,1;1,1 & |\,\,-{{\left( 2k\pi  \right)}^{2}}{{\left( -t \right)}^{\beta }}  \\
   \beta +2,\beta ,\beta ;1,1;1,1 & |\,\,-{{\left( 2k\pi  \right)}^{2}}{{\left( -t \right)}^{\beta }}  \\
\end{matrix} \right), \\
 & {{W}_{2k}}\left( t \right)={{\varphi }_{2k}}-\frac{{{\psi }_{2k}}-{{\varphi }_{2k}}}{p\,{{E}_{\beta ,2}}\left( -{{\left( 2k\pi  \right)}^{2}}{{p}^{\beta }} \right)}t\,{{E}_{\beta ,2}}\left( -{{\left( 2k\pi  \right)}^{2}}{{\left( -t \right)}^{\beta }} \right), \\
\end{aligned}\eqno (4.3)
$$
and
$$
\overline{{{\varphi }_{1k}}}={{\psi }_{1k}}-{{\varphi }_{1k}}+\frac{4k\pi {{p}^{\beta}}\left[{{\psi }_{2k}}-{{\varphi }_{2k}}\right]}{{{E}_{\beta ,2}}\left( -{{\left( 2k\pi  \right)}^{2}}{{p}^{\beta }} \right)}{{E}_{1}}\left( \begin{matrix}
   1,1;1,1 & |\,\,-{{\left( 2k\pi  \right)}^{2}}{{p}^{\beta }}  \\
   \beta +2,\beta ,\beta ;1,1;1,1 & |\,\,-{{\left( 2k\pi  \right)}^{2}}{{p}^{\beta }}  \\
\end{matrix} \right).
$$

For existence of solution, we need to prove the convergence of series corresponding to $u$, ${}_{C}D_{0t}^{\alpha }u,$ $\,{{}_{C}}D_{t0}^{\beta }u,\,\,{{u}_{xx}}$ and $f$. We will consider here the convergence of the series correspond to $u_{xx}$, since it requires stronger conditions due to the appearance of the term $\left( 2k\pi  \right)^2$. Precisely, we need to prove the convergence of the series $\sum\limits_{k=1}^{\infty }{{{\left( 2k\pi  \right)}^{2}}\left| {{V}_{ik}}\left( t \right) \right|}$  and $\sum\limits_{k=1}^{\infty }{{\left( 2k\pi  \right)}^{2}}\left| {{W}_{ik}}\left( t \right) \right|$ $\left( i=1,2 \right)$. In order to guarantee this, we impose certain conditions to given functions.

For the convergence of the series $\sum\limits_{k=1}^{\infty }{{{\left( 2k\pi  \right)}^{2}}\left| {{V}_{ik}}\left( t \right) \right|}$, we assume the following regularity conditions:
$$
\varphi \left( x \right)\in {{C}^{2}}\left[ 0,1 \right],\,\,{\varphi }'''\left( x \right)\in {{L}_{2}}\left( 0,1 \right),\,\,\varphi \left( 0 \right)=\varphi \left( 1 \right),\,\,{\varphi }'\left( 0 \right)=0.\eqno (4.4)
$$
Then on integration by parts, we get
$$
\begin{aligned}
  & \sum\limits_{k=1}^{\infty }{{{\left( 2k\pi  \right)}^{2}}\left| {{V}_{1k}}\left( t \right) \right|}=\sum\limits_{k=1}^{\infty }{{{\left( 2k\pi  \right)}^{2}}\left| {{\varphi }_{1k}} \right|}\le \sum\limits_{k=1}^{\infty }{\frac{1}{2k\pi }\left( \left| \varphi _{k}^{\left( 3 \right)} \right|+\left| \varphi _{k}^{\left( 2 \right)} \right| \right)}\le\\
  &\le \frac{1}{4}\sum\limits_{k=1}^{\infty }{\left( \frac{1}{{{\left( k\pi  \right)}^{2}}}+2{{\left| \varphi _{k}^{\left( 3 \right)} \right|}^{2}}+2{{\left| \varphi _{k}^{\left( 2 \right)} \right|}^{2}} \right)}, \\
 & \sum\limits_{k=1}^{\infty }{{{\left( 2k\pi  \right)}^{2}}\left| {{V}_{2k}}\left( t \right) \right|}=\sum\limits_{k=1}^{\infty }{{{\left( 2k\pi  \right)}^{2}}\left| {{\varphi }_{2k}} \right|}\le \sum\limits_{k=1}^{\infty }{\frac{2}{k\pi }\left( {{C}_{1}}+\left| \overline{\overline{\varphi }}_{k}^{\left( 3 \right)} \right| \right)}\le \\
 &\le \sum\limits_{k=1}^{\infty }{\left( \frac{1}{{{\left( k\pi  \right)}^{2}}}+{{C}_{2}}{{\left| \overline{\overline{\varphi }}_{k}^{\left( 3 \right)} \right|}^{2}} \right)}, \\
\end{aligned}\eqno (4.5)
$$
where $C_1,\,C_2$ are some constants and
$$
\begin{aligned}
&\varphi _{k}^{\left( 2 \right)}=2\int\limits_{0}^{1}{{\varphi }''\left( x \right)\sin 2k\pi xdx},\,\varphi _{k}^{\left( 3 \right)}=2\int\limits_{0}^{1}{{\varphi }'''\left( x \right)\left( 1-x \right)\sin 2k\pi xdx},\\
& \overline{\overline{\varphi }}_{k}^{\left( 3 \right)}=\int\limits_{0}^{1}{{\varphi }'''\left( x \right)\cos 2k\pi xdx}.\\
\end{aligned}
$$
For the convergence of the series $\sum\limits_{k=1}^{\infty }{{\left( 2k\pi  \right)}^{2}}\left| {{W}_{ik}}\left( t \right) \right|$, we first estimate ${{W}_{ik}}\left( t \right)\,\,\left( i=1,2 \right)$. Using (4.3), we have
$$
\begin{aligned}
  & \left| {{W}_{1k}}\left( t \right) \right|\le \left| {{\varphi }_{1k}} \right|+\left| \overline{{{\varphi }_{1k}}} \right|\cdot \left| t{{E}_{\beta ,2}}\left( -{{\left( 2k\pi  \right)}^{2}}{{\left( -t \right)}^{\beta }} \right) \right|+\left( \left| {{\psi }_{2k}} \right|+\left| {{\varphi }_{2k}} \right| \right)\cdot 4k\pi \times  \\
 & \times \left| \frac{{{\left( -t \right)}^{\beta +1}}{{E}_{1}}\left( \begin{matrix}
   1,1;1,1 & |\,\,-{{\left( 2k\pi  \right)}^{2}}{{\left( -t \right)}^{\beta }}  \\
   \beta +2,\beta ,\beta ;1,1;1,1 & |\,\,-{{\left( 2k\pi  \right)}^{2}}{{\left( -t \right)}^{\beta }}  \\
\end{matrix} \right)}{p{{E}_{\beta ,2}}\left( -{{\left( 2k\pi  \right)}^{2}}{{p}^{\beta }} \right)} \right|. \\
\end{aligned}\eqno (4.6)
$$
Using Theorem 2.1 and formula (2.8), one can verify that
$$
\left| {{E}_{1}}\left( \begin{matrix}
   1,1;1,1 & |\,\,-{{\left( 2k\pi  \right)}^{2}}x  \\
   {{\tau }_{1}},{{\tau }_{2}},{{\tau }_{3}};1,1;1,1 & |\,\,-{{\left( 2k\pi  \right)}^{2}}x  \\
\end{matrix} \right) \right|\le C,\eqno (4.7)
$$
for some constant $C$.

Now assuming
$$
\begin{aligned}
  & \varphi \left( x \right),\,\psi \left( x \right)\in {{C}^{3}}\left[ 0,1 \right],\,{{\varphi }^{iv}}\left( x \right),\,{{\psi }^{iv}}\left( x \right)\in {{L}_{2}}\left( 0,1 \right),\,\varphi \left( 0 \right)=\varphi \left( 1 \right),\,{\varphi }'\left( 0 \right)=0, \\
 & {\varphi }''\left( 0 \right)={\varphi }''\left( 1 \right),\psi \left( 0 \right)=\psi \left( 1 \right),\,{\psi }'\left( 0 \right)=0,\,{\psi }''\left( 0 \right)={\psi }''\left( 1 \right) \\
\end{aligned}\eqno (4.8)
$$
we obtain
$$
\begin{aligned}
  & \sum\limits_{k=1}^{\infty }{{{\left( 2k\pi  \right)}^{2}}\left| {{W}_{1k}}\left( t \right) \right|}\le \\
  &\le \sum\limits_{k=1}^{\infty }{{{C}_{3}}\left( \frac{1}{{{\left( k\pi  \right)}^{2}}}+{{\left| \varphi _{k}^{\left( 2 \right)} \right|}^{2}}+{{\left| \varphi _{k}^{\left( 3 \right)} \right|}^{2}}+{{\left| \varphi _{k}^{\left( 4 \right)} \right|}^{2}}+{{\left| \psi _{k}^{\left( 2 \right)} \right|}^{2}}+{{\left| \psi _{k}^{\left( 3 \right)} \right|}^{2}}+{{\left| \psi _{k}^{\left( 4 \right)} \right|}^{2}} \right)}, \\
 & \sum\limits_{k=1}^{\infty }{{{\left( 2k\pi  \right)}^{2}}\left| {{W}_{2k}}\left( t \right) \right|}\le \sum\limits_{k=1}^{\infty }{{{C}_{4}}\left( \frac{1}{{{\left( k\pi  \right)}^{2}}}+{{\left| \varphi _{k}^{\left( 3 \right)} \right|}^{2}}+{{\left| \psi _{k}^{\left( 3 \right)} \right|}^{2}} \right)}, \\
\end{aligned}\eqno (4.9)
$$
for some constants $C_3, \, C_4$, where
$$
\begin{aligned}
  & \varphi _{k}^{\left( 4 \right)}=\int\limits_{0}^{1}{{{\varphi }^{iv}}\left( x \right)\sin 2k\pi xdx},\,\psi _{k}^{\left( 2 \right)}=2\int\limits_{0}^{1}{{\psi }''\left( x \right)\sin 2k\pi xdx},\, \\
 & \psi _{k}^{\left( 3 \right)}=2\int\limits_{0}^{1}{{\psi }'''\left( x \right)\left( 1-x \right)\sin 2k\pi xdx},\,\psi _{k}^{\left( 4 \right)}=\int\limits_{0}^{1}{{{\psi }^{iv}}\left( x \right)\sin 2k\pi xdx}. \\
\end{aligned}
$$
Hence, the convergence of $\sum\limits_{k=1}^{\infty }{{{\left( 2k\pi  \right)}^{2}}\left| {{V}_{ik}}\left( t \right) \right|}$  and $\sum\limits_{k=1}^{\infty }{{\left( 2k\pi  \right)}^{2}}\left| {{W}_{ik}}\left( t \right) \right|$ $\left( i=1,2 \right)$ follows from the facts that $\sum\limits_{k=1}^{\infty }{\frac{1}{{{\left( k\pi  \right)}^{2}}}\,}=\frac{1}{6}$ and $\sum\limits_{k=1}^{\infty }{{{\left| g \right|}^{2}}}\le {{\left\| g \right\|}_{{{L}_{2}}\left( 0,1 \right)}}.$

Similarly, one can show the convergence of the series corresponding to $u$, $_{C}D_{0t}^{\alpha }u,\,{{\,}_{C}}D_{t0}^{\beta }u,$ and $f$. This ends the proof of the existence of the solution.

\subsection{The uniqueness}

Assume that problem (3.1)-(3.4) has two pairs of solutions $\left\{ {{u}_{1}}\left( x,t \right),\,{{f}_{1}}\left( x \right) \right\}$ and  $\left\{ {{u}_{2}}\left( x,t \right),\,{{f}_{2}}\left( x \right) \right\}$. Then the pair $\left\{u\left( x,t \right)={{u}_{1}}\left( x,t \right)-{{u}_{2}}\left( x,t \right),\,\,\,f\left( x \right)={{f}_{1}}\left( x \right)-{{f}_{2}}\left( x \right)\right\}$ satisfies the problem (3.1)-(3.4) with $\psi \left( x \right)=\varphi \left( x \right)=0$.

According to (3.35) we obtain
$$
\begin{aligned}
&{{f}_{0}}={{f}_{1k}}={{f}_{2k}}={{V}_{0}}\left( 0 \right)={{V}_{1k}}\left( 0 \right)={{V}_{2k}}\left( 0 \right)={{W}_{0}}\left( 0 \right)={{W}_{0}}^{\prime }\left( 0 \right)=\\
&={{W}_{1k}}\left( 0 \right)={{{W}'}_{1k}}\left( 0 \right)={{W}_{2k}}\left( 0 \right)={{{W}'}_{2k}}\left( 0 \right)=0.\\
\end{aligned}
$$

Based on (3.17)-(3.22), from (3.8) we deduce
$$
\begin{aligned}
  & \int\limits_{0}^{1}{}u(x,t)(1-x)dx=0,\,\int\limits_{0}^{1}{}u(x,t)(1-x)\cos 2k\pi xdx=0,\,\int\limits_{0}^{1}{}u(x,t)\sin 2k\pi xdx=0, \\
 & \int\limits_{0}^{1}{}f(x)(1-x)dx=0,\,\int\limits_{0}^{1}{}f(x)(1-x)\cos 2k\pi xdx=0,\,\int\limits_{0}^{1}{}f(x)\sin 2k\pi xdx=0,\,\,k=1,2,... \\
\end{aligned}
$$
Due to completeness of the system of functions (2.17), (2.18) in ${{L}_{2}}\left[ 0,1 \right]$, we can state that $u\left( x,t \right)=0$  a.e. in $\left[ 0,1 \right]$ for $t\in \left[ -p,q \right]$ and $f\left( x \right)=0$ a.e. in $\left[ 0,1 \right]$.

Finally, the existence and uniqueness of problem (3.1)-(3.4) can be stated as

\textbf{Theorem 4.1.} Let $0<\gamma <1.$ Suppose that the condition (4.8) holds, then the problem (3.1)-(3.4) has a unique solution set $\left\{ u\left( x,t \right),\,f\left( x \right) \right\}$ represented by (4.1)-(4.2).

\section{Special case of transmitting condition}

For the case $\gamma =1$ the condition (3.4) can be rewritten as
$$
\underset{t\to +0}{\mathop{\lim }}\,\,{{\,}_{C}}D_{0t}^{\alpha }u\left( x,t \right)=\underset{t\to -0}{\mathop{\lim }}\,\,{{u}_{t}}\left( x,t \right),\,\,\,\,\,\,0<x<1.\eqno (5.1)
$$
We note that uniqueness and existence of nontrivial solutions of the problem (3.1)-(3.3), (5.1) with $\varphi(x)=\psi(x)=0$ were studied by Salakhitdinov and Karimov [18]. Here we will present the existence of the solution for the case $\varphi(x)\ne0$ and $\psi(x)\ne 0$.

Considering (3.30)-(3.33) and using (5.1), we get
$$
\begin{aligned}
&{{f}_{0}}={{W}_{0}}^{\prime }\left( 0 \right),\,\,{{f}_{1k}}+4k\pi {{V}_{2k}}\left( 0 \right)-{{(2k\pi )}^{2}}{{V}_{1k}}\left( 0 \right)={{W}_{1k}}^{\prime }\left( 0 \right),\\
&{{f}_{2k}}-{{(2k\pi )}^{2}}{{V}_{2k}}\left( 0 \right)={{W}_{2k}}^{\prime }\left( 0 \right).\\
\end{aligned}
\eqno (5.2)
$$

From (3.23)-(3.29) and (5.2), we obtain the following two systems of algebraic equations:
$$
\left\{ \begin{aligned}
  & {{W}_{0}}(0)+\left( \frac{{{p}^{\beta }}}{\Gamma \left( \beta +1 \right)}+p \right){{W}_{0}}^{\prime }(0)={{\psi }_{0}}, \\
 & {{W}_{0}}(0)+\frac{{{q}^{\alpha }}}{\Gamma \left( \alpha +1 \right)}{{W}_{0}}^{\prime }(0)={{\varphi }_{0}}, \\
 & {{f}_{0}}={{W}_{0}}^{\prime }(0), \\
\end{aligned} \right.\eqno (5.3)
$$
and
$$
\left\{ \begin{aligned}
  & {{f}_{1k}}={{W}_{1k}}^{\prime }(0)+{{(2k\pi )}^{2}}{{W}_{1k}}(0)-4k\pi {{W}_{2k}}(0), \\
 & {{f}_{2k}}={{W}_{2k}}^{\prime }(0)+{{(2k\pi )}^{2}}{{W}_{2k}}(0), \\
 & {{W}_{2k}}(0)+{{q}^{\alpha }}{{E}_{\alpha ,\alpha +1}}\left( -{{\left( 2k\pi  \right)}^{2}}{{q}^{\alpha }} \right){{W}_{2k}}^{\prime }(0)={{\varphi }_{2k}}, \\
 & {{W}_{2k}}(0)+{{W}_{2k}}^{\prime }(0)\left[ p{{E}_{\beta ,2}}\left( -{{\left( 2k\pi  \right)}^{2}}{{p}^{\beta }} \right)+{{p}^{\beta }}{{E}_{\beta ,\beta +1}}\left( -{{\left( 2k\pi  \right)}^{2}}{{p}^{\beta }} \right) \right]={{\psi }_{2k}}, \\
 & {{W}_{1k}}\left( 0 \right)+{{{{W}'}}_{1k}}\left( 0 \right){{q}^{\alpha }}{{E}_{\alpha ,\alpha +1}}\left( -{{\left( 2k\pi  \right)}^{2}}{{q}^{\alpha }} \right)=\overline{{{\psi }_{k}}}, \\
 & {{W}_{1k}}\left( 0 \right)+{{{{W}'}}_{1k}}\left( 0 \right)\left[ p{{E}_{\beta ,2}}\left( -{{\left( 2k\pi  \right)}^{2}}{{p}^{\beta }} \right)+{{p}^{\beta }}{{E}_{\beta ,\beta +1}}\left( -{{\left( 2k\pi  \right)}^{2}}{{p}^{\beta }} \right) \right]=\widetilde{{{\psi }_{k}}}, \\
\end{aligned} \right.\eqno (5.4)
$$
where
$$
\begin{aligned}
  & \overline{{{\psi }_{k}}}={{\varphi }_{1k}}-\frac{{{\psi }_{2k}}-{{\varphi }_{2k}}}{{{\Delta }_{k}}}4k\pi {{q}^{2\alpha }}{{E}_{1}}\left( \begin{matrix}
   1,1;1,1 & |-{{(2k\pi )}^{2}}{{q}^{\alpha }}  \\
   \alpha +1,\alpha ,\alpha ;1,1;1,1 & |-{{(2k\pi )}^{2}}{{q}^{\alpha }}  \\
\end{matrix} \right), \\
 & \widetilde{{{\psi }_{k}}}={{\psi }_{1k}}-\frac{{{\psi }_{2k}}-{{\varphi }_{2k}}}{{{\Delta }_{k}}}4k\pi \,{{p}^{\beta }}\left[ p{{E}_{\beta ,2}}\left( -{{\left( 2k\pi  \right)}^{2}}{{p}^{\beta }} \right)+\right.\\
 &\left.+{{p}^{\beta }}{{E}_{1}}\left( \begin{matrix}
   1,1;1,1 & |-{{(2k\pi )}^{2}}{{p}^{\beta }}  \\
   \beta +1,\beta ,\beta ;1,1;1,1 & |-{{(2k\pi )}^{2}}{{p}^{\beta }}  \\
\end{matrix} \right) \right]. \\
\end{aligned}
$$
These systems are solvable provided the conditions
$$
\begin{aligned}
 & {{\Delta }_{k}}=p{{E}_{\beta ,2}}\left( -{{\left( 2k\pi  \right)}^{2}}{{p}^{\beta }} \right)+{{p}^{\beta }}{{E}_{\beta ,\beta +1}}\left( -{{\left( 2k\pi  \right)}^{2}}{{p}^{\beta }} \right)-{{q}^{\alpha }}{{E}_{\alpha ,\alpha +1}}\left( -{{\left( 2k\pi  \right)}^{2}}{{q}^{\alpha }} \right)\ne 0, \\
 & {{\Delta }_{0}}=p+\frac{{{p}^{\beta }}}{\Gamma \left( \beta +1 \right)}-\frac{{{q}^{\alpha }}}{\Gamma \left( \alpha +1 \right)}\ne 0. \\
\end{aligned}\eqno (5.5)
$$
Solutions of the above given systems are then given by
$$
\begin{aligned}
& {{f}_{0}}=\frac{{{\psi }_{0}}-{{\varphi }_{0}}}{{{\Delta }_{0}}}, \\
 & {{f}_{1k}}=\frac{\overline{{{\psi }_{k}}}}{{{\Delta }_{k}}}\left( {{\Delta }_{k}}-{{E}_{\alpha ,\alpha +1}}\left( -{{\left( 2k\pi  \right)}^{2}}{{q}^{\alpha }} \right)+1 \right)-\\
 &-4k\pi {{\varphi }_{2k}}-\frac{{{\psi }_{2k}}-{{\varphi }_{2k}}}{{{\Delta }_{k}}}4k\pi {{q}^{\alpha }}{{E}_{\alpha ,\alpha +1}}\left( -{{\left( 2k\pi  \right)}^{2}}{{q}^{\alpha }} \right), \\
 & {{f}_{2k}}=\frac{{{\psi }_{2k}}-{{\varphi }_{2k}}}{{{\Delta }_{k}}}+{{\left( 2k\pi  \right)}^{2}}{{\varphi }_{2k}}-\frac{{{\psi }_{2k}}-{{\varphi }_{2k}}}{{{\Delta }_{k}}}{{\left( 2k\pi  \right)}^{2}}{{q}^{\alpha }}{{E}_{\alpha ,\alpha +1}}\left( -{{\left( 2k\pi  \right)}^{2}}{{q}^{\alpha }} \right), \\
  & {{V}_{0}}\left( 0 \right)={{\varphi }_{0}}-\frac{{{q}^{\alpha }}}{\Gamma \left( \alpha +1 \right)}\frac{{{\psi }_{0}}-{{\varphi }_{0}}}{{{\Delta }_{0}}},\\
 & {{V}_{1k}}\left( 0 \right)=\frac{\overline{{{\psi }_{k}}}}{{{\Delta }_{k}}}\left( {{\Delta }_{k}}-{{E}_{\alpha ,\alpha +1}}\left( -{{\left( 2k\pi  \right)}^{2}}{{q}^{\alpha }} \right) \right),\\
  & {{V}_{2k}}\left( 0 \right)={{W}_{2k}}\left( 0 \right)={{\varphi }_{2k}}-\frac{{{\psi }_{2k}}-{{\varphi }_{2k}}}{{{\Delta }_{k}}}{{q}^{\alpha }}{{E}_{\alpha ,\alpha +1}}\left( -{{\left( 2k\pi  \right)}^{2}}{{q}^{\alpha }} \right),\\
  & {{W}_{0}}\left( 0 \right)={{\varphi }_{0}}-\frac{{{q}^{\alpha }}}{\Gamma \left( \alpha +1 \right)}\frac{{{\psi }_{0}}-{{\varphi }_{0}}}{{{\Delta }_{0}}},\\
  & {{W}_{0}}^{\prime }\left( 0 \right)=\frac{{{\psi }_{0}}-{{\varphi }_{0}}}{{{\Delta }_{0}}}, \\
  & {{W}_{1k}}\left( 0 \right)=\frac{\overline{{{\psi }_{k}}}}{{{\Delta }_{k}}}\left( {{\Delta }_{k}}-{{E}_{\alpha ,\alpha +1}}\left( -{{\left( 2k\pi  \right)}^{2}}{{q}^{\alpha }} \right) \right),\\
 & {{{{W}'}}_{1k}}\left( 0 \right)=\,\frac{\overline{{{\psi }_{k}}}}{{{\Delta }_{k}}}, \\
& {{W}_{2k}}\left( 0 \right)={{\varphi }_{2k}}-\frac{{{\psi }_{2k}}-{{\varphi }_{2k}}}{{{\Delta }_{k}}}{{q}^{\alpha }}{{E}_{\alpha ,\alpha +1}}\left( -{{\left( 2k\pi  \right)}^{2}}{{q}^{\alpha }} \right),\\
 & {{W}_{2k}}^{\prime }\left( 0 \right)=\frac{{{\psi }_{2k}}-{{\varphi }_{2k}}}{{{\Delta }_{k}}}, \\
\end{aligned}\eqno (5.7)
$$
Substituting the expressions in (5.7) into (3.17)-(3.22), we obtain
$$
\begin{aligned}
  & {{V}_{0}}\left( t \right)={{\varphi }_{0}}+\frac{{{\psi }_{0}}-{{\varphi }_{0}}}{{{\Delta }_{0}}\Gamma \left( \alpha +1 \right)}\left( {{t}^{\alpha }}-{{q}^{\alpha }} \right), \\
 & {{V}_{1k}}\left( t \right)=\left[ {{t}^{\alpha }}{{E}_{\alpha ,\alpha +1}}\left( -{{\left( 2k\pi  \right)}^{2}}{{t}^{\alpha }} \right)-{{q}^{\alpha }}{{E}_{\alpha ,\alpha +1}}\left( -{{\left( 2k\pi  \right)}^{2}}{{q}^{\alpha }} \right)+{{\Delta }_{k}} \right]\frac{\overline{{{\psi }_{k}}}}{{{\Delta }_{k}}}+ \\
 & +\frac{{{\psi }_{2k}}-{{\varphi }_{2k}}}{{{\Delta }_{k}}}4k\pi {{t}^{2\alpha }}{{E}_{1}}\left( \begin{matrix}
   1,1;1,1 & |-{{(2k\pi )}^{2}}{{t}^{\alpha }}  \\
   2\alpha +1,\alpha ,\alpha ;1,1;1,1 & |-{{(2k\pi )}^{2}}{{t}^{\alpha }}  \\
\end{matrix} \right), \\
 & {{V}_{2k}}\left( t \right)={{\varphi }_{2k}}+\frac{{{\psi }_{2k}}-{{\varphi }_{2k}}}{{{\Delta }_{k}}}\left[ {{t}^{\alpha }}{{E}_{\alpha ,\alpha +1}}\left( -{{\left( 2k\pi  \right)}^{2}}{{t}^{\alpha }} \right)-{{q}^{\alpha }}{{E}_{\alpha ,\alpha +1}}\left( -{{\left( 2k\pi  \right)}^{2}}{{q}^{\alpha }} \right) \right], \\
\end{aligned}\eqno (5.8)
$$
$$
\begin{aligned}
  & {{W}_{0}}\left( t \right)={{\varphi }_{0}}+\frac{{{\psi }_{0}}-{{\varphi }_{0}}}{{{\Delta }_{0}}}\left[ \frac{{{\left( -t \right)}^{\beta }}}{\Gamma \left( \beta +1 \right)}-t-\frac{{{q}^{\alpha }}}{\Gamma \left( \alpha +1 \right)} \right], \\
 & {{W}_{1k}}\left( t \right)=\frac{\overline{{{\psi }_{k}}}}{{{\Delta }_{k}}}\left[ {{\left( -t \right)}^{\beta }}{{E}_{\beta ,\beta +1}}\left( -{{\left( 2k\pi  \right)}^{2}}{{\left( -t \right)}^{\beta }} \right)-t{{E}_{\beta ,2}}\left( -{{\left( 2k\pi  \right)}^{2}}{{\left( -t \right)}^{\beta }} \right)-\right. \\
 & \left.-{{q}^{\alpha }}{{E}_{\alpha ,\alpha +1}}\left( -{{\left( 2k\pi  \right)}^{2}}{{q}^{\alpha }} \right)+{{\Delta }_{k}} \right]+\\
 &+\frac{{{\psi }_{2k}}-{{\varphi }_{2k}}}{{{\Delta }_{k}}}4k\pi {{\left( -t \right)}^{\beta }}\left[ {{\left( -t \right)}^{\beta }}{{E}_{1}}\left( \begin{matrix}
   1,1;1,1 & |-{{(2k\pi )}^{2}}{{\left( -t \right)}^{\beta }}  \\
   2\beta +1,\beta ,\beta ;1,1;1,1 & |-{{(2k\pi )}^{2}}{{\left( -t \right)}^{\beta }}  \\
\end{matrix} \right)- \right. \\
 & \left. -t{{E}_{1}}\left( \begin{matrix}
   1,1;1,1 & |-{{(2k\pi )}^{2}}{{\left( -t \right)}^{\beta }}  \\
   \beta +2,\beta ,\beta ;1,1;1,1 & |-{{(2k\pi )}^{2}}{{\left( -t \right)}^{\beta }}  \\
\end{matrix} \right) \right], \\
 & {{W}_{2k}}\left( t \right)={{\varphi }_{2k}}+\frac{{{\psi }_{2k}}-{{\varphi }_{2k}}}{{{\Delta }_{k}}}\left[ {{\left( -t \right)}^{\beta }}{{E}_{\beta ,\beta +1}}\left( -{{\left( 2k\pi  \right)}^{2}}{{\left( -t \right)}^{\beta }} \right)-\right. \\
 & \left. -t{{E}_{\beta ,2}}\left( -{{\left( 2k\pi  \right)}^{2}}{{\left( -t \right)}^{\beta }} \right)-{{q}^{\alpha }}{{E}_{\alpha ,\alpha +1}}\left( -{{\left( 2k\pi  \right)}^{2}}{{q}^{\alpha }} \right) \right]. \\
\end{aligned}\eqno (5.9)
$$

Thus, solution of problem (3.1)-(3.3), (5.1) is given by (3.5)-(3.7) with the coefficients given by (5.7)-(5.9). Convergence of these series can be done similarly as in the previous case ($0<\gamma<1$).

Finally, we formulate our result for this case as follows:

\textbf{Theorem 5.1.} Suppose that conditions (4.7) and (5.5) hold, then problem (3.1)-(3.3), (5.1) has a unique solution $\left\{ u\left( x,t \right),\,f\left( x \right) \right\}$, which has a representation (3.5)-(3.7) with the coefficients given by (5.7)-(5.9).

\section*{Conclusion}

In this paper, we have proved a unique solvability of an inverse-source problem for time-fractional mixed type equation with Caputo differential operator in a rectangular domain. As a main tool of investigation, we have used a series expansion of solution using bi-orthogonal system. We have shown that transmitting condition has an influence on the unique solvability of the problem. Precisely, in case of full integral form of transmitting condition ($0<\gamma <1$), solutions are obtained without any restriction on the geometry of considered domain, while in the case of the semi-integral form of transmitting condition ($\gamma =1$), there is a certain restrictions on the lower and upper bounds ($p$ and $q$) of the considered rectangle (see condition (5.5)).

During the proof of the main result, we had to simplify bulky representations and we have found a new property of recently introduced Mittag-Leffler type function as stated in Lemma 2.1.
\section*{Appendix}

\textbf{A1. Solution of (3.13)} 

According to [17, p.17] general solution of (3.13) is given by
$$
\begin{aligned}
  & {{V}_{1k}}(t)={{V}_{1k}}(0){{E}_{\alpha ,1}}\left( -{{(2k\pi )}^{2}}{{t}^{\alpha }} \right)+{{f}_{1k}}\cdot {{t}^{\alpha }}\cdot {{E}_{\alpha ,\alpha +1}}\left( -{{(2k\pi )}^{2}}{{t}^{\alpha }} \right)+ \\
 & +4k\pi \cdot {{V}_{2k}}(0)\int\limits_{0}^{t}{{{\left( t-z \right)}^{\alpha -1}}{{E}_{\alpha ,1}}\left( -{{(2k\pi )}^{2}}{{z}^{\alpha }} \right){{E}_{\alpha ,\alpha }}\left( -{{(2k\pi )}^{2}}{{\left( t-z \right)}^{\alpha }} \right)dz}+ \\
 & +4k\pi \cdot {{f}_{2k}}\int\limits_{0}^{t}{{{\left( t-z \right)}^{\alpha -1}}{{z}^{\alpha }}{{E}_{\alpha ,\alpha +1}}\left( -{{(2k\pi )}^{2}}{{z}^{\alpha }} \right){{E}_{\alpha ,\alpha }}\left( -{{(2k\pi )}^{2}}{{\left( t-z \right)}^{\alpha }} \right)dz}. \\
\end{aligned}
$$
Here we have simplified the two integrals on the right hand side of the above equation. For instance, the first integral can be simplified by setting $t-z=t\tau $ and using (2.10):
$$
\begin{aligned}
  & \int\limits_{0}^{t}{{{\left( t-z \right)}^{\alpha -1}}{{E}_{\alpha ,1}}\left( -{{(2k\pi )}^{2}}{{z}^{\alpha }} \right){{E}_{\alpha ,\alpha }}\left( -{{(2k\pi )}^{2}}{{\left( t-z \right)}^{\alpha }} \right)dz}= \\
 & ={{t}^{\alpha }}\int\limits_{0}^{1}{{{\tau }^{\alpha -1}}{{\left( 1-\tau  \right)}^{1-1}}E_{\alpha ,\alpha ;1,1}^{1,1}\left( -{{(2k\pi )}^{2}}{{t}^{\alpha }}{{\tau }^{\alpha }} \right)E_{\alpha ,1;1,1}^{1,1}\left( -{{(2k\pi )}^{2}}{{t}^{\alpha }}{{\left( 1-\tau  \right)}^{\alpha }} \right)d\tau }. \\
\end{aligned}
$$
Considering (2.2) and using (2.8) at ${{\gamma }_{1}}={{\gamma }_{2}}={{\alpha }_{1}}={{\beta }_{1}}={{\delta }_{2}}={{\delta }_{3}}={{\alpha }_{3}}={{\beta }_{3}}=1,\,{{\rho }_{1}}=\alpha ,\,{{\rho }_{2}}=1,$  $\,{{\alpha }_{2}}={{\beta }_{2}}=\alpha ,$ $x=y=-{{\left( 2k\pi  \right)}^{2}}{{t}^{\alpha }}$, we obtain 
$$
{{t}^{\alpha }}\cdot {{E}_{1}}\left( \begin{matrix}
   1,1;1,1 & |-{{(2k\pi )}^{2}}{{t}^{\alpha }}  \\
   \alpha +1,\alpha ,\alpha ;1,1;1,1 & |-{{(2k\pi )}^{2}}{{t}^{\alpha }}  \\
\end{matrix} \right).
$$
Similarly, the second integral can be simplified to get
$$
{{t}^{2\alpha }}\cdot {{E}_{1}}\left( \begin{matrix}
   1,1;1,1 & |-{{(2k\pi )}^{2}}{{t}^{\alpha }}  \\
   2\alpha +1,\alpha ,\alpha ;1,1;1,1 & |-{{(2k\pi )}^{2}}{{t}^{\alpha }}  \\
\end{matrix} \right).
$$

\textbf{A2. Integral form of solution to equation (3.14)}

According to [17, p.17] general solution of (3.14) has the following form
$$
\begin{aligned}
  & {{W}_{1k}}(t)={{W}_{1k}}(0){{E}_{\beta ,1}}\left( -{{(2k\pi )}^{2}}{{(-t)}^{\beta }} \right)-t{{W}_{1k}}^{\prime }(0){{E}_{\beta ,2}}\left( -{{(2k\pi )}^{2}}{{(-t)}^{\beta }} \right)- \\
 & -{{f}_{1k}}\int\limits_{0}^{t}{{{\left( \xi -t \right)}^{\beta -1}}{{E}_{\beta ,\beta }}\left( -{{(2k\pi )}^{2}}{{(\xi -t)}^{\beta }} \right)d\xi }- \\
 & -4k\pi \cdot {{W}_{2k}}(0)\int\limits_{0}^{t}{{{\left( \xi -t \right)}^{\beta -1}}{{E}_{\beta ,1}}\left( -{{(2k\pi )}^{2}}{{(-\xi )}^{\beta }} \right){{E}_{\beta ,\beta }}\left( -{{(2k\pi )}^{2}}{{(\xi -t)}^{\beta }} \right)d\xi }+ \\
 & +4k\pi \cdot {{W}_{2k}}^{\prime }(0)\int\limits_{0}^{t}{{{\left( \xi -t \right)}^{\beta -1}}\xi {{E}_{\beta ,2}}\left( -{{(2k\pi )}^{2}}{{(-\xi )}^{\beta }} \right){{E}_{\beta ,\beta }}\left( -{{(2k\pi )}^{2}}{{(\xi -t)}^{\beta }} \right)d\xi }- \\
 & -4k\pi \cdot {{f}_{2k}}\int\limits_{0}^{t}{{{\left( \xi -t \right)}^{\beta -1}}{{\left( -\xi  \right)}^{\beta }}{{E}_{\beta ,\beta +1}}\left( -{{(2k\pi )}^{2}}{{(-\xi )}^{\beta }} \right){{E}_{\beta ,\beta }}\left( -{{(2k\pi )}^{2}}{{(\xi -t)}^{\beta }} \right)d\xi,} \\
\end{aligned}
$$
$$
\begin{aligned}
&W_2k(t)=W_{2k}(0)E_{\beta,1}\left( -{{(2k\pi )}^{2}}{{(-t)}^{\beta }} \right)-tW{2k}'(0) E_{\beta,2}\left( -{{(2k\pi )}^{2}}{{(-t)}^{\beta }} \right)+\\
&+4k\pi f_{2k}\int\limits_0^t (\xi-t)^{\beta-1}E_{\beta,\beta}\left(-(2k\pi)^2(\xi-t)^\beta\right)\,d\xi.\\
\end{aligned}
$$

These representations can be simplified using a similar approach as in Appendix A1 to get (3.21), (3.22).

\section*{Acknowledgement}

Authors acknowledge financial support from The Research Council (TRC), Oman. This work is funded byTRC under the research agreement no. ORG/SQU/CBS/13/030.

\end{document}